\documentclass[11pt,reqno]{amsart}
\usepackage{amsmath,amsxtra,latexsym,amsthm,amssymb,amscd,pb-diagram}
\usepackage{mathrsfs,mathabx,amsfonts}
\usepackage[colorlinks=true,linkcolor=red,citecolor=blue]{hyperref}
\usepackage[margin=1in]{geometry}
\usepackage{bm}
\usepackage{color}
\usepackage{makecell,multirow,diagbox}
\usepackage{mathtools}
\usepackage[mathlines]{lineno}
\usepackage{verbatim} 
\usepackage{graphicx}
\usepackage{epsfig}
\usepackage{array}
\usepackage{enumitem}
\usepackage{silence}
\newtheorem{definition}{Definition}[section]
\newtheorem{theorem}{Theorem}[section]
\newtheorem{proposition}{Proposition}[section]
\newtheorem{lemma}{Lemma}[section]
\newtheorem{corollary}{Corollary}[section]
\newtheorem{remark}{Remark}[section]
\newtheorem{example}{Example}[section]

\begin{document} 
\title[\texorpdfstring{$\sigma$}{sigma}-evolution equations with time-dependent coefficients]{On the Cauchy problem for \texorpdfstring{$\sigma$}{sigma}-evolution equations with time-dependent coefficients via upper decay indicators}

\thanks{$^*$\textit{Corresponding author:} Tuan Anh Dao}

\maketitle
\centerline{\textbf{Duc An Phan}}
{\footnotesize
	\centerline{Department of Mathematics, Banking Academy of Vietnam}
	\centerline{12 Chua Boc, Kim Lien, Hanoi, Vietnam}
	\centerline{Email: anpd@hvnh.edu.vn}}
\medskip

\centerline{\textbf{Tuan Anh Dao}}
{\footnotesize
	\centerline{Faculty of Mathematics and Informatics, Hanoi University of Science and Technology}
	\centerline{No.1 Dai Co Viet road, Hanoi, Vietnam}
	\centerline{Email: anh.daotuan@hust.edu.vn}}

\centerline{\textbf{The Anh Cung}}
{\footnotesize
	\centerline{School of Mathematics and Computer Science, Hanoi National University of Education}
	\centerline{136 Xuan Thuy, Cau Giay, Hanoi, Vietnam}
	\centerline{Email: anhctmath@hnue.edu.vn}}

\begin{abstract}
In this paper, the first main goal is to investigate the Cauchy problem for the linear $\sigma$-evolution equation
$$
u_{tt}+a^2(t)(-\Delta)^\sigma u+b(t)u_t+m^2(t)u=0,
$$
where propagation speed, damping term and mass are time-dependent functions. Under suitable assumptions on the interaction among these coefficients, we not only derive decay estimates for the solution and its derivatives in terms of the effective parabolic time scale $\mathcal{V}_{a,b}(s,t):=\int_s^t a^2(\tau)/b(\tau) d\tau$, but also conclude the additional decay induced by the mass term and the upper decay indicators of the initial data. Especially, one can recognize that the time-derivative estimate for solutions exhibits a competition between $m^2(t)$ and $a^2(t)\left(1+\mathcal{V}_{a,b}(0,t)\right)^{-1}$. For the corresponding semi-linear problem with the nonlinear term $|u|^p$, we also establish global (in time) existence of small data solutions fulfilling decay estimates, coming from the mass factor $\exp\left(-\int_0^t m^2(\tau)/b(\tau)d\tau\right)$, under the same assumptions on these coefficients.
\end{abstract}

\tableofcontents

\section{Introduction}
\subsection{Motivation}
Let us begin with the following Cauchy problem for semi-linear $\sigma$-evolution equations with time-dependent damping:
$$
\begin{cases}
u_{t t}+\left(-\Delta\right)^\sigma u+b(t) u_t=|u|^p, & (t, x) \in[0, \infty) \times \mathbb{R}^n,\\ 
u(0, x)=u_0(x), u_t(0, x)=u_1(x), & x \in \mathbb{R}^n,
\end{cases}
$$
where $\sigma \geq 1$ and $p>1$. For the classical wave case $\sigma=1$, the role of the time-dependent damping coefficient was studied by Wirth \cite{Wirth2004, Wirth2006, Wirth2007}. His analysis distinguishes between \textit{non-effective} damping, for which the solution retains the main features of free wave propagation, and \textit{effective} damping, for which the large-time behavior becomes parabolic-like. As a continuation, effective damping conditions were considered by D'Abbicco-Lucente-Reissig \cite{DAbbicco_Reissig2013} for $\sigma=1$, and later by Aslan-Dao \cite{Dao2024} for $\sigma>1$. Afterwards, the semi-linear Cauchy problem with an effective time-dependent damping term has been extensively investigated under a variety of assumptions on the damping coefficient $b(t)$. For the model $b(t) = \mu(1+t)^{-\beta}, \mu > 0$ and $\beta \in (-1,1)$, the Fujita exponent $p_{\mathrm{Fuj}}(n):=1+2/n$ was shown in \cite{Nishihara2012} to be the critical threshold. More precisely, a blow-up result occurs when $1<p\leq p_{\mathrm{Fuj}}(n)$ and suitably nontrivial nonnegative initial data, whereas small data global solutions exist for $p_{\mathrm{Fuj}}(n)<p<\infty$ if $n=1,2$, and for $p_{\mathrm{Fuj}}(n)<p<(n+2)/(n-2)$ if $n\geq3$. Subsequently, D'Abbicco-Lucente-Reissig \cite{DAbbicco_Reissig2013} extended the global existence theory to monotone polynomial-type damping terms and small initial data in weighted energy spaces for $p_{\mathrm{Fuj}}(n)<p\le n/(n-2)$ and $n \geq 3$. Later, D'Abbicco-Lucente \cite{DAbbicco2013} developed a refined test function approach to establish blow-up results in both the subcritical regime $1 < p < p_{\mathrm{Fuj}}(n)$ and the critical case $p = p_{\mathrm{Fuj}}(n)$. For the scale-invariant case $\beta = 1$, Ikeda-Sobajima \cite{Ikeda2018} proved a blow-up result for $p_{\mathrm{Fuj}}(n) \leq p \leq p_S(n+\mu)$ with $0<\mu< (n^2+n+2)/(n+2)$, while the global existence one for the particular case $\mu = 2$ was established in \cite{DAbbicco2015, Lucente2015}. Here $p_S(n)$ denotes the positive root of $(n-1)p^2-(n+1)p-2=0$ for $n>1$, and we identify $p_S(1)=\infty$. For dispersive damping, i.e. $\beta > 1$, Lai-Takamura \cite{Takamura2018} derived blow-up results for $1 < p < p_S(n)$. On the other hand, when $\beta < -1$, Ikeda-Wakasugi \cite{Ikeda2020} showed the global existence of small data solutions for any $p > 1$ to observe the absence of a Fujita-type critical exponent in this regime.  Furthermore, Ikeda-Inui \cite{Ikeda2019} obtained sharp lifespan estimates in the critical case $p = p_{\mathrm{Fuj}}(n)$ for damping terms of the form $b(t) = (1 + t)^{-\beta}$ with $\beta \in [-1, 1)$. For more general classes of damping coefficients including $b(t) = 1 + t$ and $b(t) = (1 + t)(1 + \log(1 + t))$, Ikeda-Sobajima-Wakasugi \cite{Sobajima2019} investigated blow-up phenomena in both the critical and subcritical cases. More recently, Aslan-Dao \cite{Dao2024} studied semi-linear $\sigma$-evolution equations with effective time-dependent damping for any $\sigma\geq1$ by deriving linear decay estimates, small data global existence in the supercritical case and also catching blow-up, lifespan estimates in the subcritical and critical cases. \medskip

Besides the damping mechanism, the influence of a genuinely time-dependent propagation speed has also been investigated for classical wave equations. At the linear problem, Bui-Reissig \cite{BuiReissig2014} studied the interaction between a positive increasing propagation speed $a=a(t)$ and a time-dependent dissipation $b=b(t)$. Based on their obtained estimates, the authors in \cite{BuiSuper2015,Bui2015} considered the semi-linear Cauchy problem as follows:
$$
\begin{cases}
u_{tt}-a^2(t)\Delta u+b(t)u_t=|u|^p, & (t, x) \in [0, \infty) \times \mathbb{R}^n, \\
u(0, x)=u_0(x), \quad u_t(0, x)=u_1(x), & x \in \mathbb{R}^n,
\end{cases}
$$
where $a(t)$ is positive and increasing, $b(t)u_t$ is assumed to be effective in relation to $a(t)$. They divided the analysis into the super-exponential and sub-exponential regimes according to the behavior of the propagation speed. In the former case, global (in time) existence of small-data solutions was proved by combining Matsumura-type estimates for parameter-dependent linear Cauchy problems with a fixed-point argument. The latter case requires a different treatment of the nonlinear term, although the underlying linear estimates have the same structure, in particular, the quantitative interaction between $a(t)$ and $b(t)$ gives the range of admissible powers and the corresponding threshold exponent. Later, Aslan-Reissig \cite{AslanReissig2020} investigated semi-linear damped wave models, in which both the propagation speed and the damping coefficient may contain strong time-dependent oscillations. They showed that the unfavorable influence of such oscillations can be compensated by suitable stabilization conditions together with higher regularity assumptions on the coefficients, and obtained global (in time) existence of small data solutions in different classes of regularity. In another direction, Ebert-Marques \cite{EbertMarques2021} considered the same semi-linear structure with the decreasing propagation speed $a(t)=(1+t)^{-\ell}, \ell\in(0,1)$, and the scale-invariant damping $b(t)=\beta(1+t)^{-1}$. By exploiting an explicit representation of the corresponding linear solution through Fourier multipliers, they concluded global existence results for small data and determined the associated critical exponent. In the absence of a mass term, one recognizes that the admissible range of nonlinear powers strongly depends on the joint behavior of the propagation speed and the dissipation. \medskip

The presence of a time-dependent mass term leads to a further change in the decay mechanism. Namely, D'Abbicco-Girardi-Reissig in \cite{DAbbicco2019} studied the Cauchy problem
$$
\begin{cases}
u_{tt}-\Delta u+b(t) u_t+m^2(t) u=|u|^p, & (t, x) \in [0, \infty) \times \mathbb{R}^n, \\
u(0, x)=u_0(x), \quad u_t(0, x)=u_1(x), & x \in \mathbb{R}^n,
\end{cases}
$$
with effective time-dependent damping and a mass term satisfying $m(t)=o(b(t))$ as $t\to\infty$. With $\mathcal{V}_b(0, t):=\int_0^t b^{-1}(\tau)d\tau$ and $\chi(t):=\exp\left(-\int_0^t m^2(\tau)/b(\tau)d\tau\right)$, they derived decay estimates for the corresponding linear problem, which exhibit explicitly the influence of the mass term, moreover , the estimate for $u_t$ reflects a competition between $m^2(t)$ and $\left(1+\mathcal{V}_b(0, t)\right)^{-1}$. For the semi-linear problem, they introduced $\beta=\liminf_{t\rightarrow\infty}\mathcal{V}_b(0, t) m^2(t)$ and obtained the scale $p_\beta(n)=1+4/(n+4\beta)$ for $\beta\in[0,\infty)$, with $p_\infty(n)=1$, for small data in the energy space. Assuming additional $L^\eta$ regularity for the initial data, we can see that the corresponding scale becomes $p_{\beta,\eta}(n)=1+2\eta/(n+2\eta\beta)$. For this reason, the interaction between damping and mass affects not only the linear decay, but also the admissible range of powers in the nonlinear problem. The role of the propagation speed was established by Aslan-Ebert-Reissig \cite{AslanEbertReissig2023} by considering the following scale-invariant model:
$$
\begin{cases}
u_{tt}-(1+t)^{-2\ell}\Delta u+\dfrac{\mu_1}{1+t}u_t+\dfrac{\mu_2^2}{(1+t)^2}u=|u|^p, & (t, x) \in [0, \infty) \times \mathbb{R}^n, \\
u(0, x)=u_0(x), \quad u_t(0, x)=u_1(x), & x \in \mathbb{R}^n,
\end{cases}
$$
where $\ell>1$. The quantity $\delta=(\mu_1-1)^2-4\mu_2^2$ separates the dissipation-dominant case $\delta>0$ from the mass-dominant and balanced cases $\delta\leq0$. Their analysis provides estimates for the associated parameter-dependent linear problems and global (in time) existence results for small-data Sobolev solutions. In the dissipation-dominant regime they also obtained blow-up results and lifespan estimates to realize a critical Fujita-type exponent for suitable ranges of the parameters. More recently, Marques-Reissig \cite{MarquesReissig2026} considered a broader class of semi-linear hyperbolic models with integrable decaying speed of propagation as follows:
$$
\begin{cases}
u_{tt}-a^2(t)\Delta u+m_1\dfrac{a(t)}{\mathbb{A}(t)}u_t+m_2^2\left(\dfrac{a(t)}{\mathbb{A}(t)}\right)^2u=|u|^p, & (t, x) \in [0, \infty) \times \mathbb{R}^n, \\ u(0, x)=u_0(x), \quad u_t(0, x)=u_1(x), & x \in \mathbb{R}^n,
\end{cases}
$$
where $\mathbb{A}(t)=\int_t^\infty a(\tau)d\tau$. They proposed a classification of mass-dominant and dissipation-dominant models by means of a WKB-type analysis of the corresponding linear equations. More precisely, global (in time) well-posedness results for small-data Sobolev solutions in both the mass-dominant case and in the dissipation-dominant case with non-effective dissipation are established. From this observation, one recognizes that when all three coefficients depend on time, their effects cannot in general be separated in the large-time analysis. \medskip

Schonbek \cite{Schonbek1980, Schonbek1985, Schonbek1986} introduced the Fourier splitting method as a powerful technique for catching decay estimates for solutions to parabolic equations and the Navier-Stokes system. Later, this approach was subsequently refined by Bjorland-Schonbek in \cite{Bjorland2009} with particular emphasis on the notion of the decay character of the initial data (see Section \ref{section_Decay_characters} for the precise definition). Roughly speaking, for initial data $u_0 \in L^2$, the decay character $r^*=r^*(u_0)$ allows one to capture more accurately the large-time decay rate of solutions to the heat equation associated with the datum $u_0$. More recently, the concept of decay character has been emerged as an effective tool in the analysis of the linear problem
\begin{equation} \label{eq:linear}
\begin{cases}
u_{t t}+(-\Delta)^\sigma u+ b(t)(-\Delta)^\delta u_t=0, & (t, x) \in [0, \infty) \times \mathbb{R}^n, \\
u(0, x)=u_0(x), \quad u_t(0, x)=u_1(x), & x \in \mathbb{R}^n,
\end{cases}
\end{equation}
where $\sigma \geq 1$ and $\delta \in [0, \sigma]$, as well as its semi-linear counterpart
\begin{equation} \label{eq:semi_linear}
\begin{cases}
u_{t t}+(-\Delta)^\sigma u+ b(t)(-\Delta)^\delta u_t=|u|^p, & (t, x) \in [0, \infty) \times \mathbb{R}^n, \\
u(0, x)=u_0(x), \quad u_t(0, x)=u_1(x), & x \in \mathbb{R}^n.
\end{cases}
\end{equation}
In this direction, several contributions have appeared in recent years. For instance, Cardenas-Niche \cite{Niche2022} investigated the classical case $\sigma=1, \delta=0$ and $b(t)=1$ for \eqref{eq:linear}, Cung-Pham-Tang \cite{ADL2024} considered $\sigma \geq 1, \delta \in[0, \sigma]$ and $b(t)=1$ for \eqref{eq:linear} and \eqref{eq:semi_linear}, while Cung-Phan-Pham \cite{AAD} studied the case $\sigma \geq 1, \delta \in [0, \sigma/2]$ under effective time-dependent structural damping $b(t)(-\Delta)^\delta u_t$ for \eqref{eq:linear} and \eqref{eq:semi_linear}. These works describe decay rates of solutions to \eqref{eq:linear} in terms of the low-frequency behavior of the initial data and establish small-data global existence results for \eqref{eq:semi_linear}.

\subsection{Main purpose}
In the present paper, we would like to investigate the following Cauchy problem for linear $\sigma$-evolution equations with time-dependent coefficients:
\begin{equation}
\begin{cases} \label{eq:1.1}
u_{tt}+a^2(t)(-\Delta)^\sigma u+b(t) u_t+m^2(t) u=0, & (t, x) \in [0, \infty) \times \mathbb{R}^n, \\
u(0, x)=u_0(x), \quad u_t(0, x)=u_1(x), & x \in \mathbb{R}^n,
\end{cases}
\end{equation}
where $\sigma \geq 1$, the functions $a(t), b(t)$ and $m(t)$ stand for speed of propagation, damping and mass terms, respectively, provided that the assumptions introduced in Section \ref{Assumptions} hold. As a consequence, we also consider the corresponding semi-linear problem of \eqref{eq:1.1} in the form
\begin{equation}
\begin{cases} \label{eq:1.0}
u_{tt}+a^2(t)(-\Delta)^\sigma u+b(t) u_t+m^2(t) u=|u|^p, & (t, x) \in [0, \infty) \times \mathbb{R}^n, \\
u(0, x)=u_0(x), \quad u_t(0, x)=u_1(x), & x \in \mathbb{R}^n,
\end{cases}
\end{equation}
where $p>1$. Our first aim is to derive decay estimates for \eqref{eq:1.1} by taking into account the simultaneous effects of the time-dependent propagation speed, damping and mass (Theorem \ref{theorem:1.1}). The problem may be viewed as the $\sigma$-evolution analogue of variable-speed models mentioned in \cite{BuiReissig2014,BuiSuper2015,Bui2015}, with the additional influence of the time-dependent mass term studied in \cite{DAbbicco2019}. We carry out the linear analysis in the time-frequency phase space by using WKB-type diagonalization procedures and obtain parameter-dependent estimates for \eqref{eq:1.1} involving the quantity $\mathcal{V}_{a,b}(s,t)=\int_s^t a^2(\tau)/b(\tau)d\tau$. The mass term enters the estimates through an additional decay factor, namely, the estimate for the time derivative of solutions shows a competition between $m^2(t)$ and $a^2(t)\big(1+\mathcal{V}_{a,b}(0,t)\big)^{-1}$. Among other things, we want to describe decay estimates using decay characters in terms of considering the influence of the low-frequency behavior of the initial data. Since the decay character does not need existence for a general $L^2$ function, we work instead with the upper decay indicators. Actually, this only requires an upper bound for the Fourier mass near the origin. In this framework, our results cover some typical cases of the initial data in $L^\eta\cap L^2$ with $\eta\in[1,2]$, $L^{1,\eta}\cap L^2$ with $\eta\in[0,1]$, $\dot{H}^{-\eta}\cap L^2$ with $\eta\geq0$, as well as data in the Fourier-weighted spaces $Y^q\cap L^2$ for $q<n/2$. The second contribution of this work is to conclude the global (in time) existence of small data mild solutions to the semi-linear problem \eqref{eq:1.0} for a range of powers $p$ (Theorem \ref{theorem:1.2}) by applying parameter-dependent estimates for solutions to \eqref{eq:1.1}. The lower bound for $p$ to guarantee the global (in time) existence of small data solutions to \eqref{eq:1.0} reflects both the low-frequency behavior of the initial data and the interplay among the propagation speed, the damping and the mass term. The point worth noticing is that we do not impose any additional assumption on the coefficients for the semi-linear problem. In a comparison with previous works \cite{BuiReissig2014,BuiSuper2015,Bui2015}, our problem includes a time-dependent mass term and allows us to make a generalization of any $\sigma\geq1$, while the main one considered in \cite{DAbbicco2019} incorporates a time-dependent propagation speed and describes the low-frequency contribution through upper decay indicators without requiring the existence of decay characters. Finally, we want to say that the above-mentioned competition between the mass contribution and the effective parabolic scale in the refined estimate for $u_t$ play an important role in the nonlinear fixed-point argument. \medskip

\noindent \textbf{Notations:} Throughout the paper, we write $g\lesssim h$ when there exists a constant $C>0$ such that $g\leq Ch$, use $g\gtrsim h$ if $h\lesssim g$, and write $g\sim h$ if both $g\lesssim h$ and $g\gtrsim h$ hold. For $s \in \mathbb{R}$, the spaces $H^s$ and $\dot{H}^s$ are the nonhomogeneous and homogeneous Sobolev spaces over $\mathbb{R}^n$. We denote by $[k]^+ := \max\{k,0\}$ and by $[k]$ the integer part of $k \in \mathbb{R}$. Let us fix $\overline{r}, r_0, r_1$, given parameters from $[-n/2, \infty)$. We define $P_{r}(f)_{-}, P_{r}(f)_{+}$ as the lower  and upper decay indicators of a function $f \in L^2$ (see Definition \ref{definition_2.1}). Moreover, we denote by $r^*\left(f\right)$ the decay character associated with a function $f$ (see Definition \ref{definition_2.2}). For $r\in[-n/2, \infty)$, we introduce the quantitative upper decay quantity
$$
\Theta_r(f):=\sup_{0<\rho\leq1}\rho^{-2r-n}\int_{B(\rho)}|\widehat{f}(\xi)|^2 d \xi.
$$
If $P_r(f)_+<\infty$, then $\Theta_r(f)<\infty$, i.e. the limsup controls all sufficiently small radii, while $\|f\|_{L^2}$ controls the radii bounded away from zero. For $f\in H^s$ with $\Theta_r(f)<\infty$, we define
$$
\|f\|_{r,s} :=\|f\|_{\dot{H}^s}+\|f\|_{L^2} +\Theta_r(f)^{1/2}\,\,\text{ with } s>0, \quad \|f\|_{r,0} :=\|f\|_{L^2}+\Theta_r(f)^{1/2}.
$$
For $s_0,s_1\geq0$, we denote
$$
\mathcal{D}_{r_0, r_1}^{s_0, s_1}:=\left\{ (u_0,u_1)\in H^{s_0}\times H^{s_1}: P_{r_i}(u_i)_{+}<\infty, i=0,1 \right\},
$$
equipped with the norm $\|(u_0,u_1)\|_{\mathcal{D}_{r_0, r_1}^{s_0, s_1}}:= \|u_0\|_{r_0,s_0} +\|u_1\|_{r_1,s_1}$. Finally, we define
$$
\mathcal{V}_{a, b}(s, t) := \int_s^t \frac{a^2(\tau)}{b(\tau)} d \tau\quad \text{ for } s\in[0, t] \text{ and }\quad \chi(t):=\exp\left(-\int_0^t m^2(\tau)/b(\tau)d\tau\right).
$$

\noindent \textbf{This paper is organized as follows:} Section \ref{section_2} introduces the assumptions on the time-dependent coefficients and states the main results for the linear and semi-linear problems. Section \ref{section_3} is devoted to the proof of Theorem \ref{theorem:1.1} concerning the linear decay estimates for \eqref{eq:1.1} in terms of the upper decay indicators of the initial data. Section \ref{section_5} proves the small-data global existence result for \eqref{eq:1.0} stated in Theorem \ref{theorem:1.2}. Finally, in the Appendix, Section \ref{section_Decay_characters} recalls the definition and main properties of decay indicators and decay characters, while Section \ref{section_tools} collects some useful tools.

\section{Main results} \label{section_2}
\subsection{Assumptions} \label{Assumptions}
Let us introduce a crucial class of interactions between the time-dependent coefficients in our analysis.

\noindent\textbf{Assumption (A)}: We assume that the propagation speed $a=a(t)\in\mathcal{C}^2([0,\infty))$ satisfies the following conditions:
\begin{itemize}
\item \textbf{(A1)} $a(t)>0, a^{\prime}(t) \geq 0$ for $t \in[0, \infty)$.
\item \textbf{(A2)} $a_0 a(t)/A(t) \leq a^{\prime}(t)/a(t) \leq a_1 a(t)/A(t), a_0, a_1 \geq 0$, where $A(t)=1+\displaystyle\int_0^t a(s) d s$ is a primitive of $a(t)$.
\item \textbf{(A3)} $\left|a^{\prime \prime}(t)\right| \leq a_2 a(t)\left(a(t)/A(t)\right)^2, a_2 \geq 0$.
\end{itemize}

\noindent\textbf{Assumption (B)}: We impose the following effective-dissipation conditions on $b(t)$ with an auxiliary function $\mu=\mu(t)$:
\begin{itemize}
\item \textbf{(B1)} $b(t)>0, b(t)=\mu(t) a(t)/A(t)$.
\item \textbf{(B2)} $\mu \in \mathcal{C}^2([0, \infty))$ and $\big|\mu^{(k)}(t)\big| \leq C_k \mu(t)\left(a(t)/A(t)\right)^k$ for $k=1,2$.
\item \textbf{(B3)} $\mu(t)/A(t)$ is monotonic and $\mu(t) \rightarrow \infty$ as $t \rightarrow \infty$.
\item \textbf{(B4)} $a^2(t)/b(t)=(a(t) A(t))/\mu(t) \notin L^1\left(\mathbb{R}_+\right)$.
\item \textbf{(B5)} $a^2(t)/(A^2(t) b(t)) \in L^1\left(\mathbb{R}_+\right)$.
\end{itemize}
By Assumptions \textbf{(A2)}-\textbf{(A3)} and \textbf{(B1)}-\textbf{(B2)}, we have $b\in\mathcal{C}^2([0,\infty))$ and
$$
\left|b^{(k)}(t)\right|\lesssim b(t)\left(\frac{a(t)}{A(t)}\right)^k, \quad k=1,2.
$$

\noindent\textbf{Assumption (C)}: We assume that the non-negative mass coefficient $m=m(t)\in\mathcal{C}^2([0,\infty))$ satisfies the following conditions:
\begin{itemize}
\item \textbf{(C1)} $\left|m^{(k)}(t)\right| \leq C_k m(t)\left(a(t)/A(t)\right)^k$ for $k=1,2$. 
\item \textbf{(C2)} $m(t)=o(b(t))$ as $t \rightarrow \infty$.
\end{itemize}

We next give several examples of coefficients satisfying Assumptions \textbf{(A)}-\textbf{(B)}-\textbf{(C)}.
\begin{example}
\rm{
Let
$$
a(t)=c_0 (1+t)^{d}, \quad b(t)=\frac{c_0 c_1 (\log (c+t))^\gamma(1+t)^{d}}{1+\frac{c_0}{d+1}\left((1+t)^{d+1}-1\right)}, \quad m(t)=c_2 (1+t)^{-r},
$$
where $c_i>0$ for $i=0,1,2, d\geq0, \gamma>1, c>\exp\left(\gamma/(\min\{c_0,d+1\})\right)$ and $r\geq1$. Then, direct calculations are to claim that these coefficients satisfy Assumptions \textbf{(A1)}-\textbf{(A3)}, \textbf{(B1)}-\textbf{(B5)} and \textbf{(C1)}-\textbf{(C2)}.
}
\end{example}
\begin{example}
\rm{
Let
$$
a(t)=c_0 \mathrm{e}^t, \quad b(t)=\frac{c_0 c_1 \mathrm{e}^{(\ell+1) t}}{1+c_0\left(\mathrm{e}^t-1\right)}, \quad m(t)=c_2 (1+t)^r,
$$
where $c_i>0$ for $i=0,1,2, r>0$ and $0<\ell\leq\min\{c_0,1\}$. These coefficients satisfy Assumptions \textbf{(A1)}-\textbf{(A3)}, \textbf{(B1)}-\textbf{(B5)} and \textbf{(C1)}-\textbf{(C2)} by noticing that
$$
\frac{a^2(t)}{b(t)} \sim \mathrm{e}^{(2-\ell) t} \notin L^1\left(\mathbb{R}_+\right) \quad \text{ and } \quad \frac{a^2(t)}{A^2(t) b(t)} \sim \mathrm{e}^{-\ell t} \in L^1\left(\mathbb{R}_+\right).
$$
}
\end{example}
\begin{example}
\rm{
Let $k\geq1$. Define $A_0(t):=\mathrm{e}^t$, $A_j(t):=\mathrm{e}^{A_{j-1}(t)}$ with $j=1,\ldots,k$. Let
$$
a(t)=\prod_{j=0}^{k}A_j(t), \quad b(t)=\prod_{j=0}^{k}A_j(t), \quad m(t)=c_m(1+t)^r,
$$
where $k\geq1, c_m>0$ and $r\in\mathbb{R}$. These coefficients satisfy Assumptions \textbf{(A1)}-\textbf{(A3)}, \textbf{(B1)}-\textbf{(B5)} and \textbf{(C1)}-\textbf{(C2)}. Indeed, since $A_k^\prime(t)=a(t)$, we have
$$
A(t)=1+\int_0^t a(s) ds =1+A_k(t)-A_k(0).
$$
Since $A_k(t)\geq A_k(0)>1$ for all $t\geq0$, it follows that $A_k(t)/A_k(0) \leq A(t) \leq A_k(t)$ and
$$
\prod_{j=0}^{k-1}A_j(t) \leq \frac{a(t)}{A(t)} \leq A_k(0)\prod_{j=0}^{k-1}A_j(t).
$$
Moreover, one has $a^\prime(t)/a(t) = 1+\sum_{j=0}^{k-1}\prod_{\ell=0}^{j}A_\ell(t)$ and $1\leq \prod_{\ell=0}^{0}A_\ell(t) \leq \cdots \leq \prod_{\ell=0}^{k-1}A_\ell(t)$, so that we obtain $\prod_{j=0}^{k-1}A_j(t) \leq a^\prime(t)/a(t) \leq (k+1)\prod_{j=0}^{k-1}A_j(t)$. This means that
$$
\frac{1}{A_k(0)} \frac{a(t)}{A(t)} \leq \frac{a^\prime(t)}{a(t)} \leq (k+1)\frac{a(t)}{A(t)},
$$
which proves \textbf{(A2)}. Let $S(t):= 1+\sum_{j=0}^{k-1}\prod_{\ell=0}^{j}A_\ell(t)$. Then,
$$
a^\prime(t)=S(t)a(t) \quad \text{ and } \quad a^{\prime\prime}(t) = \left(S^\prime(t)+S^2(t)\right)a(t).
$$
Due to $S(t) \leq (k+1)\prod_{j=0}^{k-1}A_j(t)$ and
$$ S^\prime(t) \leq \frac{k(k+1)}{2} \left(\prod_{j=0}^{k-1}A_j(t)\right)^2,$$
there exists a constant $C>0$ such that
$$
|a^{\prime\prime}(t)| \leq C a(t) \left( \prod_{j=0}^{k-1}A_j(t) \right)^2 \leq CA_k^2(0) a(t) \left( \frac{a(t)}{A(t)} \right)^2.
$$
Hence, \textbf{(A3)} holds. Clearly, \textbf{(A1)} is satisfied. By choosing $\mu(t):=A(t)$ we get $\mu(t)/A(t)=1$, which is monotone, and $\mu(t)\rightarrow\infty$ as $t\to\infty$. One the one hand, one derives $\mu^\prime(t) = a(t) = \mu(t)a(t)/A(t)$ and
$$
|\mu^{\prime\prime}(t)| = |a^\prime(t)| \leq (k+1) a(t) \prod_{j=0}^{k-1}A_j(t) \leq (k+1)A_k(0) \mu(t) \left( \frac{a(t)}{A(t)} \right)^2.
$$
Thus, \textbf{(B1)}-\textbf{(B3)} hold. On the other hand, we arrive at $a^2(t)/b(t) = a(t) \notin L^1(\mathbb{R}_+)$, whereas
$$
\frac{a^2(t)}{A^2(t)b(t)} = \frac{a(t)}{A^2(t)} =-\left(\frac{1}{A(t)}\right)^\prime, \text{ i.e. } \int_0^\infty a^2(t)/(A^2(t)b(t)) dt =1,
$$
which proves \textbf{(B4)} and \textbf{(B5)}. Finally, one sees $|m^{(j)}(t)| \leq C_jm(t)(1+t)^{-j}$ with $j=1,2$ and $a(t) = \prod_{j=0}^{k}A_j(t) \geq \mathrm{e}^t$, which yields
$$
0< \frac{m(t)}{b(t)} = \frac{c_m(1+t)^r}{a(t)} \leq c_m(1+t)^r\mathrm{e}^{-t} \rightarrow 0
$$
as $t\to\infty$. Therefore, \textbf{(C1)}-\textbf{(C2)} are satisfied.
}
\end{example}

\subsection{Main theorems}
The first result concerns decay estimates for solutions to \eqref{eq:1.1}. 
\begin{theorem} \label{theorem:1.1}
Assume that $a(t), b(t)$ and $m(t)$ satisfy Assumptions \textbf{(A)}-\textbf{(B)}-\textbf{(C)}. Let $r_0,r_1\in[-n/2,\infty)$ and $(u_0,u_1)\in\mathcal{D}_{r_0, r_1}^{\alpha+j\sigma,[\alpha+(j-1)\sigma]^+}$ with $\alpha\geq0$ and $j=0,1$. Then, the problem \eqref{eq:1.1} admits a unique solution
$$
u \in \mathcal{C}\left([0, \infty), H^{\alpha+j\sigma} \right) \cap \mathcal{C}^1\left([0, \infty), H^{\alpha+(j-1) \sigma}\right),
$$
satisfying 
\begin{align} \label{estimates_linear}
\|\partial_t^j u(t,\cdot)\|_{\dot{H}^\alpha} \lesssim & \frac{\chi(t) \left(a^2(t)\left(1+\mathcal{V}_{a,b}(0,t)\right)^{-1}+m^2(t)\right)^j}{b^j(t)} \notag\\
&\times\left[\left(1+\mathcal{V}_{a,b}(0,t)\right)^{-\frac{n+2r_0+2\alpha}{4\sigma}}\|u_0\|_{r_0,\alpha+j\sigma}+\left(1+\mathcal{V}_{a,b}(0,t)\right)^{-\frac{n+2r_1+2\alpha}{4\sigma}}\|u_1\|_{r_1, [\alpha+(j-1)\sigma]^+}\right].
\end{align}
\end{theorem}
\begin{remark}
\rm{
For $j=1$, the estimate \eqref{estimates_linear} contains two competing ingredients between $m^2(t)$ and $a^2(t)\left(1+\mathcal{V}_{a,b}(0,t)\right)^{-1}$. Namely, if $m^2(t) \lesssim a^2(t)\left(1+\mathcal{V}_{a,b}(0,t)\right)^{-1}$, then
$$
\left\|u_t(t, \cdot)\right\|_{\dot{H}^\alpha} \lesssim \frac{\chi(t) a^2(t)}{b(t)}\left[\left(1+\mathcal{V}_{a, b}(0, t)\right)^{-\frac{n+ 2r_0+2\alpha}{4 \sigma}-1}\|u_0\|_{r_0, \alpha+\sigma}+\left(1+\mathcal{V}_{a, b}(0, t)\right)^{-\frac{n+ 2r_1+2\alpha}{4 \sigma}-1}\|u_1\|_{r_1, \alpha}\right].
$$
If $a^2(t)\left(1+\mathcal{V}_{a,b}(0,t)\right)^{-1}\lesssim m^2(t)$ for all sufficiently large $t$, which is the case when $\beta>0$ in \eqref{eq:inf_beta}, then the estimate \eqref{estimates_linear} follows
$$
\left\|u_t(t, \cdot)\right\|_{\dot{H}^\alpha} \lesssim \frac{\chi(t) m^2(t)}{b(t)}\left[\left(1+\mathcal{V}_{a, b}(0, t)\right)^{-\frac{n+ 2r_0+2\alpha}{4 \sigma}}\|u_0\|_{r_0, \alpha+\sigma}+\left(1+\mathcal{V}_{a, b}(0, t)\right)^{-\frac{n+ 2r_1+2\alpha}{4 \sigma}}\|u_1\|_{r_1, \alpha}\right].
$$
}
\end{remark}
\begin{remark}
\rm{
Let us give the following consequences for several standard classes of initial data from Theorem \ref{theorem:1.1}:
\begin{enumerate}[label=\textup{(\roman*)},leftmargin=2.2em]
\item If $\left(u_0, u_1\right) \in\left(H^{\alpha+j\sigma} \cap L^\eta\right) \times\left(H^{\left[\alpha+(j-1) \sigma\right]^{+}} \cap L^\eta\right)$ with $\alpha \geq 0, \eta \in [1, 2]$, then $P_{r_i}\left(u_i\right)_{+}<\infty, i=0,1$, for $r_i= -n\left(1-1/\eta\right)$ (see Example \ref{Example:2.0}). Consequently, the following estimates for $j=0,1$ are true:
$$
\begin{aligned}
\|\partial_t^j u(t,\cdot)\|_{\dot{H}^{\alpha}}\lesssim & \frac{\chi(t)\left(a^2(t) \left(1+\mathcal{V}_{a, b}(0, t)\right)^{-1}+m^2(t)\right)^j}{b^j(t)} \left(1+\mathcal{V}_{a, b}(0, t)\right)^{-\frac{n}{2\sigma}\left(\frac{1}{\eta}-\frac{1}{2}\right)-\frac{\alpha}{2\sigma}} \\
&\times \left(\|u_0\|_{H^{\alpha+j\sigma} \cap L^\eta}+\|u_1\|_{H^{\left[\alpha+(j-1) \sigma\right]^{+}}\cap L^\eta}\right).
\end{aligned}
$$
\item If $u_0, u_1$ have zero means and $\left(u_0, u_1\right) \in \left(H^{\alpha+j\sigma} \cap L^{1, \eta}\right) \times \left(H^{\left[\alpha+(j-1)\sigma\right]^+} \cap L^{1, \eta} \right)$ with $\alpha \geq 0, \eta \in [0, 1]$, then $P_{r_i}\left(u_i\right)_{+}<\infty, i=0,1$, for $r_i= \eta$ (see Example \ref{Example:2.1}). Hence, the following estimates for $j=0,1$ are true:
$$
\begin{aligned}
\|\partial_t^j u(t,\cdot)\|_{\dot{H}^{\alpha}} \lesssim & \frac{\chi(t)\left(a^2(t) \left(1+\mathcal{V}_{a, b}(0, t)\right)^{-1}+m^2(t)\right)^j}{b^j(t)} \left(1+\mathcal{V}_{a, b}(0, t)\right)^{-\frac{n+2\eta+2\alpha}{4\sigma}} \\
&\times \left(\|u_0\|_{H^{\alpha+j\sigma} \cap L^{1, \eta}}+\|u_1\|_{H^{\left[\alpha+(j-1) \sigma\right]^{+}} \cap L^{1, \eta}}\right).
\end{aligned}
$$
\item If $\left(u_0, u_1\right) \in\left(H^{\alpha+j\sigma} \cap \dot{H}^{-\eta}\right) \times\left(H^{\left[\alpha+(j-1) \sigma\right]^{+}} \cap \dot{H}^{-\eta}\right)$ with $\alpha \geq 0, \eta \geq 0$, then $P_{r_i}\left(u_i\right)_{+}<\infty, i=0,1$, for $r_i= \eta-n/2$ (see Example \ref{Example:2.2}). As a result, the following estimates for $j=0,1$ are true:
$$
\begin{aligned}
\|\partial_t^j u(t,\cdot)\|_{\dot{H}^{\alpha}}\lesssim &\frac{\chi(t)\left(a^2(t) \left(1+\mathcal{V}_{a, b}(0, t)\right)^{-1}+m^2(t)\right)^j}{b^j(t)} \left(1+\mathcal{V}_{a, b}(0, t)\right)^{-\frac{\eta + \alpha}{2\sigma}} \\
&\times \left(\|u_0\|_{H^{\alpha+j\sigma} \cap\dot{H}^{-\eta}}+\|u_1\|_{H^{\left[\alpha+(j-1) \sigma\right]^{+}}\cap\dot{H}^{-\eta}}\right).
\end{aligned}
$$
\item If $\left(u_0, u_1\right) \in\left(H^{\alpha+j\sigma}\cap Y^q\right) \times\left(H^{\left[\alpha+(j-1) \sigma\right]^{+}}\cap Y^q\right)$ with $\alpha \geq 0, q<n/2$, then $P_{r_i}\left(u_i\right)_{+}<\infty, i=0,1$, for $r_i=-q$ (see Example \ref{Example:2.3}). Therefore, the following estimates for $j=0,1$ are true:
$$
\begin{aligned}
\|\partial_t^j u(t,\cdot)\|_{\dot{H}^{\alpha}} \lesssim & \frac{\chi(t)\left(a^2(t) \left(1+\mathcal{V}_{a, b}(0, t)\right)^{-1}+m^2(t)\right)^j}{b^j(t)} \left(1+\mathcal{V}_{a, b}(0, t)\right)^{-\frac{n-2q+2\alpha}{4\sigma}} \\
&\times \left(\|u_0\|_{H^{\alpha+j\sigma} \cap Y^q}+\|u_1\|_{H^{\left[\alpha+(j-1) \sigma\right]^{+}} \cap Y^q}\right).
\end{aligned}
$$
\end{enumerate}
}
\end{remark}

The second main result is involved with global existence of small data solutions to \eqref{eq:1.0}. To state this, we introduce the following quantities:
\begin{equation}\label{eq:inf_beta}
\overline{p}(\omega,\theta) :=
\begin{cases}
    1+\left(1-\dfrac{2a_0}{2+M}\right)^{+} \dfrac{2\omega\sigma}{n+2\theta\omega\sigma} &\text{ if } \theta\in[0,\infty), \\
    1 &\text{ if }\theta=\infty,
\end{cases}\quad \beta:=\liminf_{t\to\infty}\frac{\mathcal{V}_{a,b}(0,t)m^2(t)}{a^2(t)}\in[0,\infty]
\end{equation}
for $\omega\in[1,2]$ to describe the contribution of the mass term, where $M\geq0$ is chosen such that
$$
\frac{\mu'(t)}{\mu(t)}\geq-M\frac{a(t)}{A(t)},\quad t\geq0 \quad \text{ by \textbf{(B2)}}.
$$

\begin{theorem} \label{theorem:1.2}
Assume that $a(t), b(t)$ and $m(t)$ satisfy Assumptions \textbf{(A)}-\textbf{(B)}-\textbf{(C)}. Let $\alpha \in (0, \sigma]$ and $\beta$ be as in \eqref{eq:inf_beta}. Let $r_0,r_1\in[-n/2,\infty)$ with $\min \left\{r_0, r_1 \right\} \leq 0$ and $\left(u_0,u_1\right)\in\mathcal{D}_{r_0,r_1}^{\alpha,0}$. If $p>\overline p(\omega,\beta)$ and
$$
\begin{cases}
p\in[2/\omega,\infty) &\text{ if } n\leq2\alpha,\\[1mm]
p\in[2/\omega,n/(n-2\alpha)] \quad &\text{ if } 2\alpha<n \text{ and } n(2-\omega)\leq4\alpha.
\end{cases}
$$
where $\omega:=n/(n+\min\{r_0,r_1\})$, then there exists a constant $\varepsilon>0$ such that under the assumption $\|(u_0,u_1)\|_{\mathcal{D}_{r_0,r_1}^{\alpha,0}}\leq\varepsilon$, the problem \eqref{eq:1.0} admits a unique global (in time) mild solution $u\in\mathcal{C}\left([0,\infty),H^\alpha\right)$. Moreover, the solution satisfies
$$
\begin{aligned}
\|u(t,\cdot)\|_{L^2} \lesssim &\ \chi(t) \left(1+\mathcal{V}_{a,b}(0,t)\right)^{ -\frac{n+2\min\{r_0,r_1\}}{4\sigma}} \|(u_0,u_1)\|_{\mathcal{D}_{r_0,r_1}^{\alpha,0}},\\
\|u(t,\cdot)\|_{\dot{H}^\alpha} \lesssim &\ \chi(t) \left(1+\mathcal{V}_{a,b}(0,t)\right)^{ -\frac{n+2\min\{r_0,r_1\}+2\alpha}{4\sigma}} \|(u_0,u_1)\|_{\mathcal{D}_{r_0,r_1}^{\alpha,0}}.
\end{aligned}
$$
\end{theorem}

\begin{remark}
\rm{
From Theorem \ref{theorem:1.2}, we emphasize that $\overline{p}(\omega,\beta)$ is a lower bound for the global existence result. It leaves an open problem to claim whether $\overline{p}(\omega,\beta)$ is the sharp critical exponent for \eqref{eq:1.0} or not by catching a blow-up one.
}
\end{remark}

\begin{example}
\rm{
Let us consider the constant-speed case $a(t)\equiv1$. Then, $A(t)=1+t$ and $a_0=0$, so that we obtain
$$
\overline{p}(\omega,\theta)=1+\frac{2\omega\sigma}{n+2\theta\omega\sigma} \quad \text{ as well as }\quad \overline{p}(\eta,\beta)=1+\frac{2\eta}{n+2\eta\beta}
$$
when $\sigma=1$ and $L^\eta\cap L^2$-type initial data, i.e. $\omega=\eta$. Thus, we have recovered the lower bounds appearing in \cite{DAbbicco2019}.
}
\end{example}

\begin{example}
\rm{
If $m(t)\equiv0$, then $\beta=0$ and $\chi(t)\equiv1$. For $\sigma=\alpha=1$ and $\left(u_0,u_1\right)\in \left(H^1\cap L^1\right)\times\left(L^2\cap L^1\right)$, we have $\omega=1$. We can say that our model reduces to the class of variable-speed equations considered in \cite{BuiSuper2015,Bui2015} and the lower bound for $p$ is given by
$$
\overline{p}(\omega,0)=1+\left(1-\frac{2a_0}{2+M}\right)^{+}\frac{2\omega\sigma}{n}.
$$
Hence, Theorem \ref{theorem:1.2} also allows general $\sigma\geq1$ with the initial data described by the upper decay indicators.
}
\end{example}

\begin{remark}
\rm{
The parameter $\beta$ in \eqref{eq:inf_beta} measures the contribution of the mass term to the lower bound $\overline p(\omega,\beta)$ in Theorem \ref{theorem:1.2}. The decay estimates retain the exact factor
$$
\chi(t)= \exp\left(-\int_0^t\frac{m^2(\tau)}{b(\tau)}d\tau\right).
$$
For the parameter-dependent linear problem, the corresponding mass factor is
$$
\frac{\chi(t)}{\chi(s)} = \exp\left(-\int_s^t\frac{m^2(\tau)}{b(\tau)}d\tau\right).
$$
The weighted nonlinear estimates produce a factor $\chi(s)^p$, so the Duhamel integrand contains $\chi(t)\chi(s)^{p-1}$. Proposition \ref{proposition:inf_beta} then yields the required time integrability directly from $p>\overline p(\omega,\beta)$, including the cases $\beta=0$ and $\beta=\infty$. For finite $\beta>0$, the liminf assumption \eqref{eq:inf_beta} alone is not to imply $\chi(t)\lesssim \left(1+\mathcal{V}_{a,b}(0,t)\right)^{-\beta}$. Retaining $\chi(t)$ therefore gives the decay estimates without an auxiliary power below $\beta$.
}
\end{remark}

\begin{remark}
\rm{
From Theorem \ref{theorem:1.2}, we recall two typical choices of initial data for the semi-linear problem.
\begin{enumerate}[label=\textup{(\roman*)},leftmargin=2.2em]
\item Let $\left(u_0, u_1\right) \in \left(H^{\alpha} \cap L^\eta\right) \times \left(L^2 \cap L^\eta\right)$ with $\alpha \in (0, \sigma]$ and $\eta \in [1, 2]$. Taking $r_i=-n(1-1/\eta)$, we have $P_{r_i}\left(u_i\right)_{+}<\infty, i=0,1$, and $\omega=\eta$ (see Example \ref{Example:2.0}). Hence, for sufficiently small initial data, the problem \eqref{eq:1.0} admits a unique global (in time) mild solution whenever $p>\overline{p}(\eta,\beta)$ and
$$
p \in\left[2/\eta, \infty\right), \text{ if } n \leq 2 \alpha, \quad \text{ or } \quad p \in\left[2/\eta, n/(n-2\alpha)\right], \text{ if } 2\alpha<n \text{ and } n(2-\eta)\leq4\alpha.
$$
\item Let $\left(u_0, u_1\right) \in \left(H^{\alpha} \cap \dot{H}^{-\eta}\right) \times \left(L^2 \cap \dot{H}^{-\eta}\right)$ with $\alpha \in (0, \sigma]$ and $\eta \in (0,n/2)$. Taking $r_i=\eta-n/2$, we have $P_{r_i}\left(u_i\right)_{+}<\infty, i=0,1$, and $\omega=2n/(n+2\eta)$ (see Example \ref{Example:2.2}). Hence, for sufficiently small initial data, the problem \eqref{eq:1.0} admits a unique global (in time) mild solution whenever
$$p>\overline{p}\left(\frac{2n}{n+2\eta},\beta\right)$$
and
$$
p \in\left[\frac{n+2\eta}{n}, \infty\right) \text{ if } n \leq 2 \alpha,\quad \text{ or }\quad p \in\left[\frac{n+2\eta}{n}, \frac{n}{n-2 \alpha}\right] \text{ if } 2 \alpha <n \leq \frac{\alpha(n+2\eta)}{\eta}.
$$
\end{enumerate}
}
\end{remark}

\begin{example}
\rm{
Let $a(t):=(d+1)(1+t)^d$ with $d\geq0$. Then, $A(t)=(1+t)^{d+1}$ and $a_0=d/(d+1)$. Let
$$
\mu(t):=(\log(c+t))^\gamma, \quad b(t):=\mu(t)\frac{a(t)}{A(t)}=\frac{(d+1)(\log(c+t))^\gamma}{1+t},
$$
where $\gamma>1$ and $c>1$ is sufficiently large so that Assumptions \textbf{(B1)}-\textbf{(B5)} hold. Since $\mu$ is increasing, we may take $M=0$. For any fixed $\beta>0$, choose
$$
m(t):=\sqrt{2\beta}(d+1)(1+t)^{-1}(\log(c+t))^{\frac{\gamma}{2}}.
$$
Then, Assumptions \textbf{(C1)}-\textbf{(C2)} hold and
$$
\mathcal{V}_{a, b}(0,t) =\frac{(1+t)^{2(d+1)}}{2(\log(c+t))^\gamma} \left(1+o(1)\right) \quad \text{as }t\to\infty, \quad \lim_{t\to\infty} \frac{\mathcal{V}_{a, b}(0,t)m^2(t)}{a^2(t)} =\beta.
$$
Consequently, we obtain
$$
\overline{p}(\omega,\beta)=1+\frac{2\omega\sigma}{(d+1)(n+2\beta\omega\sigma)}.
$$
Therefore, Theorem \ref{theorem:1.2} is valid whenever
$$
p>1+\frac{2\omega\sigma}{(d+1)(n+2\beta\omega\sigma)}
$$
together with the required restrictions. Thus, the factor $d+1$ reflects the influence of growing propagation speed, whereas $\beta$ quantifies a contribution of the mass term to the lower bound for the global existence range.
}
\end{example}
\begin{example}
\rm{
Let $a(t):=\mathrm{e}^t$, so that $A(t)=\mathrm{e}^t$ and $a_0=1$. Choose $\mu(t):=\mathrm{e}^{\kappa t}$, $b(t):=\mu(t)=\mathrm{e}^{\kappa t}$ with $0<\kappa\leq1$, and let $m(t):=c_m(1+t)^r$ with $c_m>0$ and $r\in\mathbb{R}$. Then, Assumptions \textbf{(B1)}-\textbf{(B5)} and \textbf{(C1)}-\textbf{(C2)} are satisfied. Since $\mu$ is increasing, we may take $M=0$, while
$$
\mathcal{V}_{a, b}(0,t)=\frac{\mathrm{e}^{(2-\kappa)t}-1}{2-\kappa}, \quad \beta=\liminf _{t\rightarrow\infty}\frac{\mathcal{V}_{a, b}(0,t)m^2(t)}{a^2(t)}=0.
$$
Consequently, $\overline{p}(\omega,0)=1$ and Theorem \ref{theorem:1.2} is valid for every $p>1$ satisfying the required restrictions.
}
\end{example}

\section{Estimates for the linear equation}
\label{section_3}
Throughout this section, for $B\in\mathbb{C}^{2\times2}$, we define $\|B\|:=\displaystyle\sup_{z\in\mathbb{C}^2\setminus\{0\}} |Bz|/|z|,
$ where $|\cdot|$ denotes the Euclidean norm on $\mathbb{C}^2$. Let us consider the following parameter-dependent Cauchy problem for $s\ge 0$:
\begin{equation}
\begin{cases} \label{eq:3.1}
u_{t t}+a^2(t) (-\Delta)^\sigma u+b(t) u_t+m^2(t) u=0, & (t, x) \in[s, \infty) \times \mathbb{R}^n, \\
u(s, x)=u_0(s,x), \quad u_t(s, x)=u_1(s,x), & (s, x) \in[0, \infty) \times \mathbb{R}^n.
\end{cases}
\end{equation}

\begin{proposition} \label{proposition:5.3}
Assume that $a(t), b(t)$ and $m(t)$ satisfy Assumptions \textbf{(A)}-\textbf{(B)}-\textbf{(C)}. Let $0\leq\alpha\leq\sigma$ and $j=0,1$. If $u_0(s, \cdot)=0$ and $u_1(s, \cdot) \in H^{\left[\alpha+(j-1) \sigma\right]^{+}}$, then the problem \eqref{eq:3.1} admits a unique solution
$$
u \in \mathcal{C}\left([s, \infty), H^{\alpha+j\sigma} \right) \cap \mathcal{C}^1\left([s, \infty), H^{\left[\alpha+(j-1) \sigma\right]^{+}}\right)
$$
satisfying the following estimates for any $t\ge s$:
\begin{align} \label{estimates_linear_u_0=0_L2}
\|\partial_t^j u(t,\cdot)\|_{\dot{H}^\alpha} \lesssim &\frac{\chi(t)}{\chi(s)b(s)} \left(\frac{m^2(t)}{b(t)}+b(t)\left(1+\frac{b^2(t)}{a^2(t)} \mathcal{V}_{a, b}(s, t)\right)^{-1}\right)^j \notag\\
&\times \left(\left(\frac{A(s)}{\mu(s)}\right)^2+\mathcal{V}_{a, b}(s, t)\right)^{-\frac{(1-j)\alpha}{2 \sigma}} \|u_1(s, \cdot)\|_{H^{\left[\alpha+(j-1) \sigma\right]^{+}}}.
\end{align}
In addition, if $P_{r_1}\left(u_1(s,\cdot)\right)_{+}<\infty$ for some $r_1\in[-n/2,\infty)$, then for all $t\geq s$ we have
\begin{align} \label{estimates_linear_u_0=0}
\|\partial_t^j u(t,\cdot)\|_{\dot{H}^\alpha} &\lesssim \frac{\chi(t)}{\chi(s)b(s)} \left(\frac{m^2(t)}{b(t)}+b(t)\left(1+\frac{b^2(t)}{a^2(t)} \mathcal{V}_{a, b}(s, t)\right)^{-1}\right)^j \notag\\
&\quad\times \Bigg\{\left(1+\mathcal{V}_{a, b}(s, t)\right)^{-\frac{n+2r_1+2j\alpha}{4 \sigma}}\left(\left(\frac{A(s)}{\mu(s)}\right)^2+\mathcal{V}_{a, b}(s, t)\right)^{-\frac{(1-j)\alpha}{2\sigma}} \notag\\
&\qquad +\left(\frac{\mu(s)}{A(s)}\right)^{\frac{(1-j)\alpha}{\sigma}}\exp\left(-c\int_s^t b(\tau)d\tau\right)\Bigg\}\|u_1(s, \cdot)\|_{r_1, \left[\alpha+(j-1) \sigma\right]^{+}},
\end{align}
where $c>0$ denotes a sufficiently small constant independent of $s$ and $t$.
\end{proposition}

If $u_0(s, \cdot)=0, u_1 \in H^{\left[\alpha+(j-1) \sigma\right]^{+}} \cap L^\eta$ with $\alpha \geq 0$ and $\eta \in [1, 2]$, then $P_{r_1}\left(u_1\right)_{+}<\infty$ for $r_1= -n(1-1/\eta)$ (see Example \ref{Example:2.0}). As a consequence of Proposition \ref{proposition:5.3}, we obtain the following result.

\begin{corollary} \label{corollary:linear}
Assume that $a(t), b(t)$ and $m(t)$ satisfy Assumptions \textbf{(A)}-\textbf{(B)}-\textbf{(C)}. Let $0\leq\alpha\leq\sigma$ and $j=0,1$. If $u_0(s, \cdot)=0$ and $u_1(s, \cdot) \in H^{\left[\alpha+(j-1) \sigma\right]^{+}} \cap L^\eta$ with $\eta \in [1, 2]$, then the Sobolev solution to \eqref{eq:3.1} satisfies the following estimates for all $t\geq s$:
$$
\begin{aligned}
\|\partial_t^j u(t,\cdot)\|_{\dot{H}^\alpha} &\lesssim \frac{\chi(t)}{\chi(s)b(s)} \left(\frac{m^2(t)}{b(t)}+b(t)\left(1+\frac{b^2(t)}{a^2(t)} \mathcal{V}_{a, b}(s, t)\right)^{-1}\right)^j \\
&\quad \times \Bigg\{\left(1+\mathcal{V}_{a, b}(s, t)\right)^{-\frac{n}{2\sigma}\left(\frac{1}{\eta}-\frac{1}{2}\right)-\frac{j\alpha}{2\sigma}}\left(\left(\frac{A(s)}{\mu(s)}\right)^2+\mathcal{V}_{a, b}(s, t)\right)^{-\frac{(1-j)\alpha}{2\sigma}} \\
&\qquad +\left(\frac{\mu(s)}{A(s)}\right)^{\frac{(1-j)\alpha}{\sigma}}\exp\left(-c\int_s^t b(\tau)d\tau\right)\Bigg\}\|u_1(s, \cdot)\|_{H^{\left[\alpha+(j-1) \sigma\right]^{+}}\cap L^\eta}
\end{aligned}
$$
and the following $L^2$-$L^2$ estimates for $t\ge s$:
$$
\begin{aligned}
\|\partial_t^j u(t,\cdot)\|_{\dot{H}^\alpha} \lesssim &\frac{\chi(t)}{\chi(s)b(s)} \left(\frac{m^2(t)}{b(t)}+b(t)\left(1+\frac{b^2(t)}{a^2(t)} \mathcal{V}_{a, b}(s, t)\right)^{-1}\right)^j \\
&\times \left(\left(\frac{A(s)}{\mu(s)}\right)^2+\mathcal{V}_{a, b}(s, t)\right)^{-\frac{(1-j)\alpha}{2\sigma}}\|u_1(s, \cdot)\|_{H^{\left[\alpha+(j-1) \sigma\right]^{+}}}.
\end{aligned}
$$
\end{corollary}
After applying the Fourier transform with respect to $x$, we obtain from \eqref{eq:3.1} the resulting equation
\begin{equation}
\begin{cases} \label{eq:3.2}
\widehat{u}_{t t}+a^2(t)|\xi|^{2 \sigma} \widehat{u}+b(t) \widehat{u}_t + m^2(t)\widehat{u}=0, & (t, \xi) \in[s, \infty) \times \mathbb{R}^n, \\
\widehat{u}(s, \xi)= \widehat{u}_0(s,\xi), \quad \widehat{u}_t(s, \xi)=\widehat{u}_1(s,\xi), & (s, \xi) \in[0, \infty) \times \mathbb{R}^n.
\end{cases}
\end{equation}
Let us now introduce
$$
\widehat{v}(t, \xi)=\frac{\nu(t)}{\nu(s)} \widehat{u}(t, \xi) \quad \text{ with } \nu(t)=\exp \left(\frac{1}{2} \int_0^t b(\tau) d \tau\right),
$$
so that \eqref{eq:3.2} can be expressed in the form
\begin{equation} \label{eq:3.3}
\begin{cases}
\widehat{v}_{t t}+\mathcal{M}(t, \xi) \widehat{v}=0, & (t, \xi) \in[s, \infty) \times \mathbb{R}^n, \\
\widehat{v}(s, \xi)= \widehat{u}_0(s, \xi), \quad \widehat{v}_t(s, \xi)= b(s)\widehat{u}_0(s, \xi)/2 +\widehat{u}_1(s, \xi), & (s, \xi) \in[0, \infty) \times \mathbb{R}^n,
\end{cases}
\end{equation}
where
\begin{equation} \label{M(t,xi)}
\mathcal{M}(t, \xi):=a^2(t)|\xi|^{2 \sigma}- b^2(t)/4-b^{\prime}(t)/2 + m^2(t).
\end{equation}
By \textbf{(B1)}, we may rewrite \eqref{M(t,xi)} as follows:
$$
\mathcal{M}(t, \xi)=a^2(t)|\xi|^{2\sigma}-\frac{1}{4} \mu^2(t) \frac{a^2(t)}{A^2(t)}-\frac{1}{2}\left(\mu(t) \frac{a(t)}{A(t)}\right)^{\prime}+ m^2(t).
$$
Assumptions \textbf{(A2)}, \textbf{(B2)} and \textbf{(B3)} show that $b^{\prime}(t)$ is a negligible term in comparison with $b^2(t)$, i.e. $\left|b^{\prime}(t)\right|=o\left(b^2(t)\right)$ as $t \rightarrow \infty$. Indeed, we have
\begin{align} \label{ineq:to_0}
\frac{\left|b^{\prime}(t)\right|}{b^2(t)} =\frac{\left|\left(\mu(t) \frac{a(t)}{A(t)}\right)^{\prime}\right|}{\mu^2(t) \frac{a^2(t)}{A^2(t)}} &=\frac{\left|\mu^{\prime}(t) \frac{a(t)}{A(t)}+\mu(t) \frac{a^{\prime}(t) A(t)-a^2(t)}{A^2(t)}\right|}{\mu^2(t) \frac{a^2(t)}{A^2(t)}} \leq \frac{\left|\mu^{\prime}(t)\right| \frac{a(t)}{A(t)}+C_1 \mu(t) \frac{a^2(t)}{A^2(t)}}{\mu^2(t) \frac{a^2(t)}{A^2(t)}} \notag\\
&\le \frac{C \mu(t) \frac{a^2(t)}{A^2(t)}+C_1 \mu(t) \frac{a^2(t)}{A^2(t)}}{\mu^2(t) \frac{a^2(t)}{A^2(t)}} \lesssim \frac{1}{\mu(t)} \rightarrow 0 \text { for } t \rightarrow \infty.
\end{align}
By \eqref{ineq:to_0} and \textbf{(C2)}, we have $\left|b^{\prime}(t)\right|+m^2(t)=o\left(b^2(t)\right)$ as $t\rightarrow\infty$. Hence, $a^2(t)|\xi|^{2\sigma}-b^2(t)/4$ is the principal part of $\mathcal{M}(t,\xi)$. We define the separating line
$$
\Gamma:=\left\{(t, \xi) \in[0, \infty) \times \mathbb{R}^n: a(t)|\xi|^{\sigma}=b(t)/2\right\},
$$
which divides the extended phase space into two regions as follows:
$$
\begin{aligned}
\Pi_{\text {hyp }}:=&\left\{(t, \xi) \in[0, \infty) \times \mathbb{R}^n: a(t)|\xi|^{\sigma}> b(t)/2 \right\}, \\ 
\Pi_{\mathrm{ell}}:=&\left\{(t, \xi) \in[0, \infty) \times \mathbb{R}^n: a(t)|\xi|^{\sigma}< b(t)/2\right\}.
\end{aligned}
$$
Let us define the auxiliary weight function
$$
\langle\xi\rangle_{a,b}(t):=\sqrt{\left| a^2(t) |\xi|^{2 \sigma}-b^2(t)/4 \right|}
$$
to split the extended phase space into the following zones:
$$
\begin{aligned}
Z_{\mathrm{hyp}}:=&\left\{(t, \xi) \in\{t \geq 0\} \times \mathbb{R}^n:\langle\xi\rangle_{a,b}(t) \geq \epsilon b(t)/2\right\} \cap \Pi_{\text {hyp }}, \\
Z_{\mathrm{ell}}:=&\left\{(t, \xi) \in\{t \geq 0\} \times \mathbb{R}^n:\langle\xi\rangle_{a,b}(t) \geq \epsilon b(t)/2\right\} \cap \Pi_{\mathrm{ell}}, \\
Z_{\mathrm{red}}:=&\left\{(t, \xi) \in\{t \geq 0\} \times \mathbb{R}^n:\langle\xi\rangle_{a,b}(t) \leq \epsilon b(t)/2\right\},
\end{aligned}
$$
where $\epsilon>0$ will be a sufficiently small constant. Hence, it holds that $\langle\xi\rangle_{a,b}(t) \sim a(t) |\xi|^{\sigma}$ in $Z_{\mathrm{hyp}}$ and $\langle\xi\rangle_{a,b}(t) \sim b(t)$ in $Z_{\mathrm{ell}}$. The diagonalization procedure below is restricted to $t\geq s\geq t_0$, with $t_0>0$ sufficiently large. The remaining finite-time cases follow from standard energy estimates together with the composition property of fundamental solutions.

\subsection{Considerations in the hyperbolic zone}
In $Z_{\mathrm{hyp}}$ we have
$\langle\xi\rangle_{a,b}(t)\sim a(t)|\xi|^\sigma$. We introduce the micro-energy $V=\left(\langle\xi\rangle_{a,b}(t) \widehat{v}, \mathcal{D}_t \widehat{v}\right)^{\mathrm{T}}$, which satisfies
\begin{equation} \label{eq:3.5}
\mathcal{D}_t V=\left(\begin{array}{cc} 0
& \langle\xi\rangle_{a,b}(t) \\
\langle\xi\rangle_{a,b}(t) & 0
\end{array}\right) V +\left(\begin{array}{cc}
\frac{\mathcal{D}_t\langle\xi\rangle_{a,b}(t)}{\langle\xi\rangle_{a,b}(t)} & 0 \\
\frac{-b^{\prime}(t)+2 m^2(t)}{2\langle\xi\rangle_{a,b}(t)} & 0
\end{array}\right) V.
\end{equation}
Hereafter, we use the notation $\mathcal{D}_t=-i \partial_t$.
\begin{lemma} \label{lemma3.1}
For any sufficiently small $\epsilon>0$, the fundamental solution $\mathcal{E}_{\mathrm{hyp}}=\mathcal{E}_{\mathrm{hyp}}(t, s, \xi)$ to \eqref{eq:3.5} satisfies
$$
\|\mathcal{E}_{\mathrm{hyp}}(t, s, \xi)\| \lesssim \left(\frac{\nu(t)}{\nu(s)}\right)^{2\epsilon},
$$
for $t\geq s\geq t_0$ and $(t, \xi),(s, \xi) \in Z_{\mathrm{hyp}}$.
\end{lemma}
\begin{proof}
We diagonalize the principal part as in \cite[Lemma 3.3]{DAbbicco2019}. By \textbf{(A2)}, \textbf{(B1)}-\textbf{(B3)}, \textbf{(C2)} and \eqref{ineq:to_0}, the remainder $R(t, \xi)$ satisfies
$$
\|R(t, \xi)\| \lesssim \frac{a(t)}{A(t)}+\frac{|b^{\prime}(t)|+m^2(t)}{\langle\xi\rangle_{a,b}(t)}=o(b(t))
$$
uniformly in $Z_{\mathrm{hyp}}$. For sufficiently large $t_0$, we obtain the desired estimate by Gronwall's inequality.
\end{proof}

\subsection{Considerations in the elliptic zone}
In this zone, we introduce the auxiliary function
$$
\psi(t,\xi):=\sqrt{b^2(t)/4-a^2(t)|\xi|^{2\sigma}-m^2(t)}.
$$
Indeed, from the definition of $Z_{\mathrm{ell}}$ and Assumption \textbf{(C2)}, for all sufficiently large $t_0$ we have
$$
a^2(t)|\xi|^{2\sigma}+m^2(t)\leq (1- \epsilon^2/2)b^2(t)/4,
$$
so it holds that $\psi(t,\xi)\sim b(t)$ uniformly in $Z_{\mathrm{ell}}$. Moreover, \eqref{M(t,xi)} can be written as $\mathcal{M}(t,\xi)=-\psi^2(t,\xi)- b'(t)/2$. We take the micro-energy $V=\left(\psi(t,\xi)\widehat{v},\mathcal{D}_t\widehat{v}\right)^T$ to derive
\begin{equation} \label{eq:3.6}
\mathcal{D}_tV=\psi(t,\xi)\left(\begin{array}{cc}
0&1\\
-1&0
\end{array}\right)V+
\left(\begin{array}{cc}
\frac{\mathcal{D}_t\psi(t,\xi)}{\psi(t,\xi)}&0\\
-\frac{b'(t)}{2\psi(t,\xi)}&0
\end{array}\right)V.
\end{equation}

\begin{lemma} \label{lemma3.2}
For $t\geq s\geq t_0$, the fundamental solution $\mathcal{E}_{\mathrm{ell}}=\mathcal{E}_{\mathrm{ell}}(t,s,\xi)$ to \eqref{eq:3.6} satisfies
$$
\left\|\mathcal{E}_{\mathrm{ell}}(t,s,\xi)\right\|
\lesssim
\exp\left(\int_s^t\left(\psi(\tau,\xi)+\frac{\partial_\tau\psi(\tau,\xi)}{2\psi(\tau,\xi)}+\frac{b'(\tau)}{4\psi(\tau,\xi)}\right)d\tau\right)
$$
for $(t,\xi),(s,\xi)\in Z_{\mathrm{ell}}$.
\end{lemma}

\begin{proof}
We follow the diagonalization argument in
\cite[Section 3.3.2]{BuiReissig2014}. We introduce the matrices
$$
Q=\left(
\begin{array}{cc}
i&-i\\ 1&1
\end{array}\right), \quad Q^{-1}=\frac{1}{2}\left(
\begin{array}{cc}
-i&1\\ i&1
\end{array}\right),
$$
and carry out the change of variable $V^{(0)}:=Q^{-1}V$. Then, one has $\mathcal{D}_tV^{(0)}=[D(t,\xi)+R(t,\xi)]V^{(0)}$, where
$$
D(t,\xi)=i\psi(t,\xi)\left(\begin{array}{cc} -1&0\\ 0&1
\end{array}\right), \qquad
R(t,\xi)=\frac{\mathcal{D}_t\psi(t,\xi)}{2\psi(t,\xi)}\left(\begin{array}{cc}
1&-1\\ -1&1
\end{array}\right)
-\frac{i b'(t)}{4\psi(t,\xi)}
\left(\begin{array}{cc}
1&-1\\ 1&-1
\end{array}\right).
$$
By \textbf{(A2)}-\textbf{(A3)}, \textbf{(B1)}-\textbf{(B2)} and \textbf{(C1)}, together with $\psi(t,\xi)\sim b(t)$, we obtain
$$
\left|\partial_t^k\psi(t,\xi)\right| \lesssim \psi(t,\xi) \left(\frac{a(t)}{A(t)}\right)^k, \quad k=1,2.
$$
Consequently, we have
$$
\|R(t,\xi)\| \lesssim \frac{a(t)}{A(t)}, \quad \|\mathcal{D}_t R(t,\xi)\| \lesssim \left(\frac{a(t)}{A(t)}\right)^2.
$$
We perform one more diagonalization to obtain an integrable remainder by denoting $F^{(0)}(t,\xi):=\operatorname{diag}R(t,\xi)$ and choosing the off-diagonal matrix $N^{(1)}=N^{(1)}(t,\xi)$ such that
$$
[D(t,\xi),N^{(1)}(t,\xi)]
+R(t,\xi)-F^{(0)}(t,\xi)=0.
$$
We put $d_1(t,\xi):=-i\psi(t,\xi)$ and $d_2(t,\xi):=i\psi(t,\xi)$. Then we have
$$
N^{(1)}_{jk}(t,\xi)=-\frac{R_{jk}(t,\xi)}{d_j(t,\xi)-d_k(t,\xi)},
\quad j\neq k.
$$
Since $|d_1(t,\xi)-d_2(t,\xi)|=2\psi(t,\xi)\sim b(t)$, we obtain $\|N^{(1)}(t,\xi)\| \lesssim a(t)/(A(t)b(t))=1/\mu(t)$. Hence $N_1(t,\xi):=\mathrm{I}_2+N^{(1)}(t,\xi)$ is invertible for all sufficiently large $t_0$, and both $N_1$ and $N_1^{-1}$ are uniformly bounded. Let us now introduce $V^{(1)}:=N_1^{-1}V^{(0)}$. Thus, one has $\mathcal{D}_tV^{(1)}=\left(D+F^{(0)}+R_1\right)V^{(1)}$, where
$$
R_1=N_1^{-1} \left(RN^{(1)}-N^{(1)}F^{(0)}-\mathcal{D}_tN^{(1)}\right).
$$
Because $\|N^{(1)}(t,\xi)\| \lesssim a(t)/(A(t)b(t))=1/\mu(t)$, we have
$$
\|R(t,\xi)N^{(1)}(t,\xi)\|+\|N^{(1)}(t,\xi)F^{(0)}(t,\xi)\|\lesssim \frac{a(t)}{A(t)\mu(t)}=\frac{a^2(t)}{A^2(t)b(t)}.
$$
On the other hand, differentiating
$$
N^{(1)}_{jk}(t,\xi)=-\frac{R_{jk}(t,\xi)}
{d_j(t,\xi)-d_k(t,\xi)}, \quad j\neq k,
$$
and using $|d_1(t,\xi)-d_2(t,\xi)|=2\psi(t,\xi)\sim b(t)$, together with
$$
\|\mathcal{D}_tR(t,\xi)\| \lesssim \left(\frac{a(t)}{A(t)}\right)^2, \quad |\partial_t\psi(t,\xi)| \lesssim b(t)\frac{a(t)}{A(t)},
$$
we achieve $\|\mathcal{D}_tN^{(1)}(t,\xi)\|
\lesssim a^2(t)/\left(A^2(t)b(t)\right)$. Hence, it entails $\|R_1(t,\xi)\| \lesssim a^2(t)/\left(A^2(t)b(t)\right)$. For $t\geq s\geq t_0$ with $(\tau,\xi)\in Z_{\mathrm{ell}}$ for all $\tau\in[s,t]$, from \textbf{(B5)}, we obtain
$$
\int_s^t\|R_1(\tau,\xi)\| d\tau \lesssim \int_s^t\frac{a^2(\tau)}{A^2(\tau)b(\tau)} d\tau \leq \int_{t_0}^{\infty} \frac{a^2(\tau)}{A^2(\tau)b(\tau)} d\tau<\infty
$$
uniformly in $s$, $t$ and $\xi$. For $t\geq s\geq t_0$, we denote by $\mathcal{E}_{\mathrm{ell},1}=\mathcal{E}_{\mathrm{ell},1}(t,s,\xi)$ the fundamental solution to the transformed system after the second diagonalization, i.e.
$$
\mathcal{E}_{\mathrm{ell},1}(t,s,\xi)=\exp\left(\int_s^t\left(\psi(\tau,\xi)+\frac{\partial_\tau\psi(\tau,\xi)}{2\psi(\tau,\xi)}+\frac{b'(\tau)}{4\psi(\tau,\xi)}\right)d\tau\right)
Q_{\mathrm{ell},1}(t,s,\xi)
$$
for $(t,\xi),(s,\xi)\in Z_{\mathrm{ell}}$, where $Q_{\mathrm{ell},1}(t,s,\xi)$ is a uniformly bounded matrix. The two diagonal terms generated by $D+F^{(0)}$ correspond to
$$
\psi(t,\xi)+\frac{\partial_t\psi(t,\xi)}{2\psi(t,\xi)}+\frac{b'(t)}{4\psi(t,\xi)} \quad \text{ and } \quad
-\psi(t,\xi)+\frac{\partial_t\psi(t,\xi)}{2\psi(t,\xi)}-\frac{b'(t)}{4\psi(t,\xi)},
$$
which lead to
$$
-2\psi(t,\xi)-\frac{b'(t)}{2\psi(t,\xi)}=-2\psi(t,\xi)\left(1+\frac{b'(t)}{4\psi^2(t,\xi)}\right).
$$
Since $\psi(t,\xi)\sim b(t)$ and $|b'(t)|=o(b^2(t))$, for sufficiently large $t_0$ we have
$$
-2\psi(t,\xi)-\frac{b'(t)}{2\psi(t,\xi)}<0,\quad t\geq t_0.
$$
Applying Gronwall's inequality to the corresponding Volterra integral equation, we obtain
$$
\|Q_{\mathrm{ell},1}(t,s,\xi)\|  \lesssim \exp\left(\int_s^t\|R_1(\tau,\xi)\|d\tau\right)\lesssim1.
$$
Consequently,
$$
\begin{aligned}
\|\mathcal{E}_{\mathrm{ell}}(t,s,\xi)\|=&\|QN_1(t,\xi)\mathcal{E}_{\mathrm{ell},1}(t,s,\xi)N_1^{-1}(s,\xi)Q^{-1}\| \lesssim \|\mathcal{E}_{\mathrm{ell},1}(t,s,\xi)\|\\
\lesssim&\exp\left(\int_s^t\left(\psi(\tau,\xi)+\frac{\partial_\tau\psi(\tau,\xi)}{2\psi(\tau,\xi)}+\frac{b'(\tau)}{4\psi(\tau,\xi)}\right)d\tau\right).
\end{aligned}
$$
Hence, the proof of Lemma \ref{lemma3.2} is complete.
\end{proof}

\subsection{Considerations in the reduced zone}
In $Z_{\mathrm{red}}$, we take the micro-energy $V=\left(\epsilon b(t) \widehat{v}, \mathcal{D}_t \widehat{v}\right)^{\mathrm{T}}$. By \eqref{eq:3.3}, $V$ satisfies
\begin{equation} \label{eq:3.9}
\mathcal{D}_t V=\left(\begin{array}{cc}
\frac{\mathcal{D}_t b(t)}{b(t)} & \epsilon b(t) \\
\frac{\mathcal{M}(t, \xi)}{\epsilon b(t)} & 0
\end{array}\right) V.
\end{equation}
\begin{lemma} \label{lemma3.5}
The fundamental solution $\mathcal{E}_{\mathrm{red}}=\mathcal{E}_{\mathrm{red}}(t, s, \xi)$ to \eqref{eq:3.9} satisfies the following estimate:
$$
\left\|\mathcal{E}_{\mathrm{red}}(t, s, \xi)\right\| \lesssim\left(\frac{\nu(t)}{\nu(s)}\right)^{4 \epsilon},
$$
for $t \geq s \geq t_0$ with a sufficiently large $t_0$ and $(t, \xi),(s, \xi) \in Z_{\mathrm{red}}$.
\end{lemma}
\begin{proof}
We follow the argument in \cite[Lemma 3.4]{DAbbicco2019}. By \eqref{ineq:to_0}, Assumption \textbf{(C2)} and the definition of $Z_{\mathrm{red}}$, taking $0<\epsilon<1/4$ we have
$$
|b'(t)|+m^2(t)\leq \epsilon^2b^2(t)/4, \quad \langle\xi\rangle_{a,b}(t)\leq \epsilon b(t)/2, \quad t\geq t_0.
$$
Since $\big|a^2(t)|\xi|^{2\sigma}-b^2(t)/4\big|=\langle\xi\rangle^2_{a,b}(t)$, it follows that
\begin{equation}\label{eq:red-M-bound}
|\mathcal{M}(t,\xi)|\leq\langle\xi\rangle^2_{a,b}(t)+ |b'(t)|/2+m^2(t)\leq \epsilon^2 b^2(t)/2.
\end{equation}
Using $|b'(t)|/b(t)\leq\epsilon^2 b(t)/4$ and \eqref{eq:red-M-bound}, we obtain
$$
\left\|
\begin{pmatrix}
\frac{\mathcal{D}_t b(t)}{b(t)} & \epsilon b(t)\\
\frac{\mathcal{M}(t,\xi)}{\epsilon b(t)} & 0
\end{pmatrix}
\right\| \leq \frac{|b'(t)|}{b(t)}+\epsilon b(t) \leq 2\epsilon b(t).
$$
Therefore, the application of Gronwall's inequality gives the desired estimate.
\end{proof}

\subsection{Proof of Theorem \ref{theorem:1.1} and Proposition \ref{proposition:5.3}}

\medskip
\noindent \textbf{Estimates for the fundamental solution.} Let us now introduce
$$
d(t, s, \xi)=\frac{\nu(s)}{\nu(t)} \frac{h(s, \xi)}{h(t, \xi)}\|\mathcal{E}(t, s, \xi)\|.
$$
Here $h(t,\xi)$ is chosen according to the zone. If $(t,\xi)$ and $(s,\xi)$ belong to the same zone $Z_j$ with $j\in\{\mathrm{hyp},\mathrm{ell},\mathrm{red}\}$, we write $d_j(t,s,\xi)$. Moreover, we have
$$
h(t, \xi) \sim \begin{cases} a(t) |\xi|^{\sigma} & (t, \xi) \notin Z_{\mathrm{ell}}, \\ b(t) & (t, \xi) \in Z_{\mathrm{ell}}.\end{cases}
$$
Lemmas \ref{lemma3.1} and \ref{lemma3.5} give
$$
d_{\mathrm{hyp}}(t, s, \xi) \lesssim \frac{a(s)}{a(t)} \left(\frac{\nu(s)}{\nu(t)}\right)^{1-2 \epsilon}, \quad d_{\mathrm{red}}(t, s, \xi) \lesssim \frac{a(s)}{a(t)} \left(\frac{\nu(s)}{\nu(t)}\right)^{1-4 \epsilon}.
$$
For the elliptic zone, we have the following estimate.
\begin{lemma} \label{lemma3.6}
For $t\geq s\geq t_0$ and $(t,\xi),(s,\xi)\in Z_{\mathrm{ell}}$, we have
$$
d_{\mathrm{ell}}(t, s, \xi) \lesssim \exp \left(-|\xi|^{2 \sigma} \mathcal{V}_{a, b}(s, t)-\int_s^t \frac{m^2(\tau)}{b(\tau)} d \tau\right).
$$
\end{lemma}
\begin{proof}
An analogous estimate is proved in \cite[Lemma 3.7]{DAbbicco2019}. Using $\psi^2(t,\xi)=b^2(t)/4-a^2(t)|\xi|^{2\sigma}-m^2(t)$ and the coefficient bounds, we obtain
$$
\begin{aligned}
\psi(t,\xi)-\frac{b(t)}{2}+\frac{b'(t)-2\partial_t\psi(t,\xi)}{4\psi(t,\xi)} \leq &-\frac{a^2(t)|\xi|^{2\sigma}+m^2(t)}{b(t)}-\frac{\left(a^2(t)|\xi|^{2\sigma}+m^2(t)\right)^2}{b^3(t)} \\
&+C\frac{a(t)}{A(t)}\frac{a^2(t)|\xi|^{2\sigma}+m^2(t)}{b^2(t)} \\
\leq & -\frac{a^2(t)|\xi|^{2\sigma}+m^2(t)}{b(t)}+C\frac{a^2(t)}{A^2(t)b(t)}.
\end{aligned}
$$
Here we used $\psi(t,\xi)\sim b(t)$ and Young's inequality. Integrating this inequality and using Lemma \ref{lemma3.2}, the definition of $d_{\mathrm{ell}}$ and \textbf{(B5)}, we obtain the desired estimate.
\end{proof}
The previous estimates for the fundamental matrix give
$$
\begin{aligned}
|\widehat{u}(t, \xi)| \lesssim & d(t, s, \xi)\left(|\widehat{u}_0(s, \xi)|+\frac{1}{h(s, \xi)}|\widehat{u}_1(s, \xi)|\right), \\
\left|\widehat{u}_t(t, \xi)\right| \lesssim & (h(t, \xi)+b(t)) d(t, s, \xi)\left(|\widehat{u}_0(s, \xi)|+\frac{1}{h(s, \xi)}|\widehat{u}_1(s, \xi)|\right).
\end{aligned}
$$
For the second estimate, we used
$$
\frac{\nu(t)}{\nu(s)}\widehat{u}_t=\widehat{v}_t-\frac{b(t)}{2}\widehat{v}.
$$

\medskip
\noindent \textbf{Estimates for the solution.}
Since the estimates in $Z_{\mathrm{hyp}}$ and $Z_{\mathrm{red}}$ have the same form, we define
$$
\Lambda(t)=\left(\frac{\sqrt{1-\epsilon^2}\mu(t)}{2A(t)}\right)^{\frac{1}{\sigma}},
$$
so that $(t,\xi)\in Z_{\mathrm{ell}} \Longleftrightarrow |\xi|\leq\Lambda(t)$. By \textbf{(B3)}, $\Lambda(t)$ is monotonic. For $t\geq s\geq0$, we define
$$
\Lambda_0(t, s)=\min \{\Lambda(t), \Lambda(s)\}, \quad \Lambda_\infty(t,s)=\max\{\Lambda(t),\Lambda(s)\}.
$$
We now turn to the low-frequency contribution.
\begin{lemma} \label{lemma3.7}
Let $w(s, \cdot) \in L^2$ and $P_{\overline{r}}\left(w\right)_{+}<\infty$ with $\overline{r} \in [-n/2, \infty)$. Then
\begin{equation} \label{eq:3.10}
\left\||\xi|^\alpha \mathrm{e}^{-|\xi|^{2 \sigma} \mathcal{V}_{a, b}(s, t)} \widehat{w}(s, \cdot)\right\|_{L^2\left(|\xi| \leq \Lambda_0(t, s)\right)} \lesssim \left(1+\mathcal{V}_{a, b}(s, t)\right)^{-\frac{n+ 2 \overline{r}+2\alpha}{4 \sigma}} \|w(s, \cdot)\|_{\overline{r}, 0},
\end{equation}
for any $\alpha\geq0$ and $\mathcal{V}_{a,b}(s,t)\geq1$. The same estimate holds for all $0\leq s\leq t$ if either $s=0$ or $\mu(t)/A(t)$ is non-increasing.
\end{lemma}
\begin{proof}
For $\mathcal{V}_{a,b}(s,t)\geq1$, we have
$$
|\xi|^\alpha \mathrm{e}^{-|\xi|^{2\sigma}\mathcal{V}_{a,b}(s,t)} \lesssim \mathcal{V}_{a,b}(s,t)^{-\frac{\alpha}{2\sigma}} \mathrm{e}^{-\frac{1}{2}|\xi|^{2\sigma}\mathcal{V}_{a,b}(s,t)}.
$$
By Plancherel's theorem and Proposition \ref{decay_1}, we obtain \eqref{eq:3.10}. If $0\leq \mathcal{V}_{a,b}(s,t)\leq1$, then $\Lambda_0(t,0)\leq\Lambda(0)$ for $s=0$, while the same bound holds for arbitrary $s$ when $\mu(t)/A(t)$ is non-increasing. Hence,
$$
\left\||\xi|^\alpha \mathrm{e}^{-|\xi|^{2\sigma}\mathcal{V}_{a,b}(s,t)}\widehat{w}(s,\cdot)\right\|_{L^2(|\xi|\leq\Lambda_0(t,s))}\leq \Lambda(0)^\alpha\|w(s,\cdot)\|_{L^2}\lesssim(1+\mathcal{V}_{a,b}(s,t))^{-\frac{n+2\overline{r}+2\alpha}{4\sigma}}\|w(s,\cdot)\|_{\overline{r},0}.
$$
This completes the proof.
\end{proof}

We now combine the low-, intermediate- and high-frequency estimates. For the initial velocity, we keep the factor $h^{-1}(s,\xi)$. For later use, we first note that, for any fixed $N>0$ and $c>0$,
\begin{equation} \label{eq:exp-parabolic}
\exp\left(-c\int_s^t b(\tau)d\tau\right)\leq C_N \left(1+\left(\frac{\mu(s)}{A(s)}\right)^2\mathcal{V}_{a,b}(s,t)\right)^{-N}.
\end{equation}
Indeed, by \textbf{(B2)} and the definition of $M$, we have
$$
\mu(\tau)\geq\mu(s)\left(\frac{A(\tau)}{A(s)}\right)^{-M}, \quad s\leq\tau\leq t,
$$
which entails
$$
\left(\frac{\mu(s)}{A(s)}\right)^2\mathcal{V}_{a,b}(s,t)=\mu^2(s)\int_s^t\frac{a(\tau)}{A(\tau)}
\left(\frac{A(\tau)}{A(s)}\right)^2\frac{1}{\mu(\tau)}d\tau \leq\frac{\mu(s)}{2+M}
\left(\left(\frac{A(t)}{A(s)}\right)^{2+M}-1\right).
$$
If $\log(A(t)/A(s))\leq1$, Assumption \textbf{(B2)} gives $\mu(\tau)\sim\mu(s)$ on $[s,t]$, and hence
$$
\left(\frac{\mu(s)}{A(s)}\right)^2\mathcal{V}_{a,b}(s,t)
\lesssim
\mu(s)\log\left(\frac{A(t)}{A(s)}\right)
\lesssim
\int_s^t b(\tau)d\tau.
$$
Let now $\log(A(t)/A(s))\geq1$. The same local comparison on the part of the interval for which $A(\tau)/A(s)\leq\mathrm{e}$ yields
$$
\int_s^t b(\tau)d\tau\gtrsim\mu(s),
$$
whereas
$$
\int_s^t b(\tau)d\tau
\geq
\left(\inf_{\theta\geq s}\mu(\theta)\right)\log\left(\frac{A(t)}{A(s)}\right).
$$
It is obvious to see that $\displaystyle\lim_{s\rightarrow\infty}\inf_{\theta\geq s}\mu(\theta)= \infty$ since $\mu(t)\rightarrow\infty$. Consequently, for any $\delta>0$, we obtain
$$
\log\left(1+\left(\frac{\mu(s)}{A(s)}\right)^2\mathcal{V}_{a,b}(s,t)\right) \leq \delta\int_s^t b(\tau)d\tau+C_\delta.
$$
When $\log(A(t)/A(s))\geq1$ and $s$ is sufficiently large, the previous estimates yield
$$
\log(1+\mu(s)) \lesssim \int_s^t b(\tau) d\tau, \quad \log\frac{A(t)}{A(s)} \lesssim \int_s^t b(\tau) d\tau,
$$
by the previous two lower bounds. For bounded $s$, we split the integral at a fixed large time and absorb the compact part into the constant. Taking $\delta>0$ sufficiently small, depending on $c$ and $N$, gives \eqref{eq:exp-parabolic}. No monotonicity of $\mu$ is needed.

For $0\leq\alpha\leq\sigma$, which is sufficient for Proposition
\ref{proposition:5.3}, we have in the elliptic zone $h(s,\xi)\sim b(s)$ and $|\xi|\leq\Lambda_0(t,s)\leq\Lambda(s) \sim (\mu(s)/A(s))^{1/\sigma}$. Hence,
$$
|\xi|^\alpha \mathrm{e}^{-|\xi|^{2\sigma}\mathcal{V}_{a,b}(s,t)}
\lesssim
\left(\left(\frac{A(s)}{\mu(s)}\right)^2+\mathcal{V}_{a,b}(s,t)\right)^{-\frac{\alpha}{2\sigma}},
$$
and Lemma \ref{lemma3.6} yields
\begin{align} \label{eq:3.11}
\left\||\xi|^\alpha d(t,s,\xi)\frac{1}{h(s,\xi)}\widehat{w}(s,\cdot)\right\|_{L^2(|\xi|\leq\Lambda_0(t,s))} &\lesssim \frac{\chi(t)}{\chi(s)b(s)}
\left(\left(\frac{A(s)}{\mu(s)}\right)^2+\mathcal{V}_{a,b}(s,t)\right)^{-\frac{\alpha}{2\sigma}}
\|w(s,\cdot)\|_{L^2}.
\end{align}
By \textbf{(C2)}, for any $\delta>0$ we have $m^2(t)/b(t)\leq\delta b(t)$ for $t\geq t_0$, provided that $t_0$ is sufficiently large. Hence, the factor $\chi(t)/\chi(s)$ can be included in the exponential estimates in the hyperbolic and reduced zones. We denote by $c_\epsilon>0$ the remaining exponential constant. For high frequencies, $|\xi|\geq\Lambda_\infty(t,s)$ and $h(s,\xi)\sim a(s)|\xi|^\sigma$. Hence,
$$
\frac{|\xi|^\alpha}{h(s,\xi)}
\lesssim\frac{1}{b(s)}\left(\frac{\mu(s)}{A(s)}\right)^{\frac{\alpha}{\sigma}}.
$$
Combining this with the hyperbolic and reduced estimates and using \eqref{eq:exp-parabolic}, we obtain \eqref{eq:3.11}.

It remains to consider the intermediate frequencies. If $\mu(t)/A(t)$ is non-decreasing, then $\Lambda(s)\leq|\xi|\leq\Lambda(t)$. Let $t_\xi\in[s,t]$ be defined by $\Lambda(t_\xi)=|\xi|$. The trajectory first passes through the hyperbolic/reduced part and then through the elliptic part. Consequently, it yields
$$
\begin{aligned}
|\xi|^\alpha d(t,s,\xi)\frac{1}{h(s,\xi)}
\lesssim&\frac{\chi(t)}{\chi(s)b(s)}
|\xi|^\alpha
\exp\left(-c_\epsilon\int_s^{t_\xi}b(\tau)d\tau
-|\xi|^{2\sigma}\mathcal{V}_{a,b}(t_\xi,t)\right).
\end{aligned}
$$
If $\mu(t)/A(t)$ is non-increasing, then $\Lambda(t)\leq|\xi|\leq\Lambda(s)$, and the order of the zones is reversed. In this case, it holds
$$
\begin{aligned}
|\xi|^\alpha d(t,s,\xi)\frac{1}{h(s,\xi)}
\lesssim&\frac{\chi(t)}{\chi(s)b(s)}
|\xi|^\alpha
\exp\left(-|\xi|^{2\sigma}\mathcal{V}_{a,b}(s,t_\xi)
-c_\epsilon\int_{t_\xi}^t b(\tau)d\tau\right).
\end{aligned}
$$
We show that both cases give the same bound. If $\mu(t)/A(t)$ is non-decreasing, writing $L_1=\log(A(t_\xi)/A(s))$, Assumption \textbf{(B2)} gives a constant $C>0$ such that
$$
\frac{\mu(t_\xi)}{A(t_\xi)}\frac{A(s)}{\mu(s)}\leq\mathrm{e}^{CL_1}.
$$
Since $\mu(t)\rightarrow\infty$ and $b(t)=\mu(t)a(t)/A(t)$, we choose $t_0$ sufficiently large so that
$$
C\frac{\alpha}{\sigma}L_1\leq\frac{c_\epsilon}{2}\int_s^{t_\xi}b(\tau)d\tau,
$$
for $s\geq t_0$. Hence, using $|\xi|^\sigma\sim \mu(t_\xi)/A(t_\xi)$ and then \eqref{eq:exp-parabolic} one has
$$
\begin{aligned}
|\xi|^\alpha\exp\left(-c_\epsilon\int_s^{t_\xi}b(\tau)d\tau\right)
&\lesssim \left(\frac{\mu(s)}{A(s)}\right)^{\frac{\alpha}{\sigma}}\exp\left(-\frac{c_\epsilon}{2}\int_s^{t_\xi}b(\tau)d\tau\right)\\
&\lesssim \left(\frac{\mu(s)}{A(s)}\right)^{\frac{\alpha}{\sigma}}\left(1+\left(\frac{\mu(s)}{A(s)}\right)^2\mathcal{V}_{a,b}(s,t_\xi)\right)^{-\frac{\alpha}{2\sigma}}
=\left(\left(\frac{A(s)}{\mu(s)}\right)^2+\mathcal{V}_{a,b}(s,t_\xi)\right)^{-\frac{\alpha}{2\sigma}}.
\end{aligned}
$$
The remaining compact-time part is absorbed into the constant by the composition property of fundamental solutions. On the other hand, we know
$$
|\xi|^\alpha\mathrm{e}^{-|\xi|^{2\sigma}\mathcal{V}_{a,b}(t_\xi,t)}\lesssim \mathcal{V}_{a,b}(t_\xi,t)^{-\frac{\alpha}{2\sigma}}.
$$
Taking the minimum of the last two bounds gives
$$
|\xi|^\alpha\exp\left(-c_\epsilon\int_s^{t_\xi}b(\tau)d\tau-|\xi|^{2\sigma}\mathcal{V}_{a,b}(t_\xi,t)\right)
\lesssim\left(\left(\frac{A(s)}{\mu(s)}\right)^2+\mathcal{V}_{a,b}(s,t_\xi)+\mathcal{V}_{a,b}(t_\xi,t)\right)^{-\frac{\alpha}{2\sigma}}.
$$
If $\mu(t)/A(t)$ is non-increasing, the elliptic part gives
$$
|\xi|^\alpha\mathrm{e}^{-|\xi|^{2\sigma}\mathcal{V}_{a,b}(s,t_\xi)}\lesssim\left(\left(\frac{A(s)}{\mu(s)}\right)^2+\mathcal{V}_{a,b}(s,t_\xi)\right)^{-\frac{\alpha}{2\sigma}},
$$
while \eqref{eq:exp-parabolic}, applied on $[t_\xi,t]$, yields
$$
|\xi|^\alpha\exp\left(-c_\epsilon\int_{t_\xi}^t b(\tau)d\tau\right)\lesssim\left(\left(\frac{A(t_\xi)}{\mu(t_\xi)}\right)^2+\mathcal{V}_{a,b}(t_\xi,t)\right)^{-\frac{\alpha}{2\sigma}}.
$$
Since $\mu(t_\xi)/A(t_\xi)\leq \mu(s)/A(s)$, the maximum of $\left(A(s)/\mu(s)\right)^2+\mathcal{V}_{a,b}(s,t_\xi)$ and $\left(A(t_\xi)/\mu(t_\xi)\right)^2+\mathcal{V}_{a,b}(t_\xi,t)$ controls $\left(A(s)/\mu(s)\right)^2+\mathcal{V}_{a,b}(s,t_\xi)+\mathcal{V}_{a,b}(t_\xi,t)$. Hence the same bound follows. Therefore, in both cases we arrive at
\begin{equation} \label{eq:3.15_L2}
\left\||\xi|^\alpha d(t,s,\xi)\frac{1}{h(s,\xi)}\widehat{w}(s,\cdot)\right\|_{L^2}
\lesssim\frac{\chi(t)}{\chi(s)b(s)}
\left(\left(\frac{A(s)}{\mu(s)}\right)^2+\mathcal{V}_{a,b}(s,t)\right)^{-\frac{\alpha}{2\sigma}}
\|w(s,\cdot)\|_{L^2},
\end{equation}
which holds uniformly for $0\leq s\leq t$. In the elliptic zone, Proposition \ref{decay_1} gives
$$
\left\||\xi|^\alpha \mathrm{e}^{-|\xi|^{2\sigma}\mathcal{V}_{a,b}(s,t)}\widehat{w}(s,\cdot)\right\|_{L^2(|\xi|\leq\Lambda_0(t,s))}
\lesssim
\left(1+\mathcal{V}_{a,b}(s,t)\right)^{-\frac{n+2\overline{r}}{4\sigma}}
\left(\left(\frac{A(s)}{\mu(s)}\right)^2+\mathcal{V}_{a,b}(s,t)\right)^{-\frac{\alpha}{2\sigma}}
\|w(s,\cdot)\|_{\overline{r},0}.
$$
The estimates in the hyperbolic and reduced zones decay exponentially. Combining the preceding estimates and using \textbf{(B2)} and $\mu(t)\to\infty$ as in \eqref{eq:3.15_L2}, we obtain
\begin{align} \label{eq:3.15}
\left\||\xi|^\alpha d(t,s,\xi)\frac{1}{h(s,\xi)}\widehat{w}(s,\cdot)\right\|_{L^2} &\lesssim \frac{\chi(t)}{\chi(s)b(s)}
\left\{
\left(1+\mathcal{V}_{a,b}(s,t)\right)^{-\frac{n+2\overline{r}}{4\sigma}}
\left(\left(\frac{A(s)}{\mu(s)}\right)^2+\mathcal{V}_{a,b}(s,t)\right)^{-\frac{\alpha}{2\sigma}} \right. \notag\\
&\hspace{2cm} \left. + \left(\frac{\mu(s)}{A(s)}\right)^{\frac{\alpha}{\sigma}}
\exp\left(-c\int_s^t b(\tau)d\tau\right)
\right\}
\|w(s,\cdot)\|_{\overline{r},0}.
\end{align}
Here and below $c>0$ denotes a sufficiently small constant independent of $s$ and $t$. For $s=0$, \eqref{eq:exp-parabolic} shows that the exponential term is dominated by any fixed power of $1+\mathcal{V}_{a,b}(0,t)$, while $\left(A(0)/\mu(0)\right)^2$ is fixed. Thus, \eqref{eq:3.15} reduces to
$$
\left\||\xi|^\alpha d(t,0,\xi)\frac{1}{h(0,\xi)}\widehat{w}(0,\cdot)\right\|_{L^2}
\lesssim\frac{\chi(t)}{\chi(0)b(0)}
\left(1+\mathcal{V}_{a,b}(0,t)\right)^{-\frac{n+2\overline{r}+2\alpha}{4\sigma}}
\|w(0,\cdot)\|_{\overline{r},0}.
$$

\medskip
\noindent \textbf{Estimates for the time derivative.}
Let us next refine the estimate for $\|\partial_tu(t,\cdot)\|_{\dot{H}^\alpha}$ when $\mathcal{V}_{a,b}(s,t)$ is large. Recalling
$$
u(t, x)=\varphi_0(t, s, x) *_{(x)} u_0(s, x)+\varphi_1(t, s, x) *_{(x)} u_1(s, x),
$$
we can observe that the multiplier $\widehat{\varphi}_i(t,s,\xi)$ satisfies
$$
\widehat{\varphi}_i^{\prime \prime}+\left(a^2(t) |\xi|^{2 \sigma}+m^2(t)\right) \widehat{\varphi}_i+b(t) \widehat{\varphi}_i^{\prime}=0, \quad \widehat{\varphi}_i(s, s, \xi)=\delta_i^0, \quad \widehat{\varphi}_i^{\prime}(s, s, \xi)=\delta_i^1.
$$
We define $\Psi_i(t,s,\xi):=\widehat{\varphi}_i'(t,s,\xi)$. Then, it holds
$$
\Psi_i^{\prime}+b(t) \Psi_i=-\left(a^2(t) |\xi|^{2 \sigma}+m^2(t)\right) \widehat{\varphi}_i, \quad \Psi_i(s, s, \xi)=\delta_i^1.
$$
Hence, one arrives at
$$
\Psi_i(t, s, \xi)=\frac{\nu^2(s)}{\nu^2(t)} \delta_i^1-\int_s^t \frac{\nu^2(\tau)}{\nu^2(t)}\left(a^2(\tau) |\xi|^{2 \sigma}+m^2(\tau)\right) \widehat{\varphi}_i(\tau, s, \xi) d \tau.
$$
Let $(t,\xi),(s,\xi)\in Z_{\mathrm{ell}}$. The previous elliptic estimate gives
$$
\left|\Psi_i(t, s, \xi)\right| \lesssim \frac{\nu^2(s)}{\nu^2(t)} \delta_i^1+I,
$$
where
$$
I:=\int_s^t \frac{\nu^2(\tau)}{\nu^2(t)}\left(\frac{a^2(\tau) |\xi|^{2 \sigma}+m^2(\tau)}{b^i(s)}\right) \exp \left(-\int_s^\tau \frac{a^2(\theta) |\xi|^{2 \sigma}+m^2(\theta)}{b(\theta)} d \theta\right) d \tau.
$$
We follow the integration-by-parts argument leading to (54) in \cite{DAbbicco2019}. Put $Q(\tau,\xi):=a^2(\tau)|\xi|^{2\sigma}+m^2(\tau)$. By \textbf{(A1)}, \textbf{(B1)}-\textbf{(B3)}, \textbf{(C1)} and the estimate for $b'$, we have
$$
-\partial_\tau\left(\frac{Q(\tau,\xi)}{b(\tau)}\right)+\frac{Q^2(\tau,\xi)}{b^2(\tau)}\leq \left(\frac{C}{\mu(\tau)}+\frac14\right)Q(\tau,\xi)\leq \frac12Q(\tau,\xi)
$$
on $[s,t]$ for sufficiently large $t_0$. Here we used $Q(\tau,\xi)\leq b^2(\tau)/4$ in $Z_{\mathrm{ell}}$. Integrating by parts in $I$ and discarding the non-positive lower endpoint term, we obtain
$$
I\leq \frac{Q(t,\xi)}{b(t)b^i(s)} \exp\left(-\int_s^t\frac{Q(\tau,\xi)}{b(\tau)}d\tau\right)+\frac{I}{2}.
$$
Therefore, we have shown that
\begin{equation} \label{eq:3.16}
\left|\Psi_i(t, s, \xi)\right| \lesssim \frac{\nu^2(s)}{\nu^2(t)} \delta_i^1+\frac{a^2(t)|\xi|^{2 \sigma}+m^2(t)}{b(t) b^i(s)} \exp \left(-|\xi|^{2 \sigma} \mathcal{V}_{a, b}(s, t)-\int_s^t \frac{m^2(\tau)}{b(\tau)} d \tau\right).
\end{equation}
By \textbf{(C2)}, we have $m(t)\lesssim b(t)$ for $t\geq t_0$ with $t_0$ sufficiently large. The estimate on bounded time intervals follows by adjusting the implicit constant. Because $a(t)|\xi|^\sigma\lesssim b(t)$ in $Z_{\mathrm{ell}}$, the estimate \eqref{eq:3.16} improves the previous bound for $u_t$ and
$$
\left|\widehat{u}_t(t, \xi)\right| \lesssim\left|\Psi_0(t, s, \xi)\right||\widehat{u}_0(s, \xi)|+\left|\Psi_1(t, s, \xi)\right||\widehat{u}_1(s, \xi)|,
$$
for $(t,\xi),(s,\xi)\in Z_{\mathrm{ell}}$. We shall use the following variant of Lemma \ref{lemma3.7}.
\begin{lemma} \label{lemma3.8}
Let $w(s, \cdot) \in L^2$ and $ P_{\overline{r}}\left(w\right)_{+}<\infty$ with $\overline{r} \in \left[-n/2, \infty\right)$. Then we have
\begin{align} \label{eq:3.17}
&\left\|\frac{|\xi|^\alpha\left(a^2(t) |\xi|^{2 \sigma}+m^2(t)\right)}{b(t) b^j(s)} \mathrm{e}^{-|\xi|^{2 \sigma} \mathcal{V}_{a, b}(s, t)} \widehat{w}(s, \cdot)\right\|_{L^2\left(|\xi| \leq \Lambda_0(t, s)\right)} \notag\\
&\qquad \lesssim \frac{a^2(t) \mathcal{V}_{a, b}^{-1}(s, t)+m^2(t)}{b(t) b^j(s)}\left(1+\mathcal{V}_{a, b}(s, t)\right)^{-\frac{n+ 2 \overline{r}+2\alpha}{4 \sigma}}\|w(s, \cdot)\|_{\overline{r}, 0},
\end{align}
for any $\alpha \geq 0, s \geq 0, t>s$ with $\mathcal{V}_{a,b}(s,t)\geq1$, and $j=0, 1$.
\end{lemma}
\begin{proof}
Proceeding as in the proof of Lemma \ref{lemma3.7}, we obtain the following estimates:
$$
\begin{aligned}
\left\||\xi|^{\alpha+2\sigma} \mathrm{e}^{-|\xi|^{2 \sigma} \mathcal{V}_{a, b}(s, t)} \widehat{w}(s, \cdot)\right\|_{L^2\left(|\xi| \leq \Lambda_0(t, s)\right)} \lesssim & \left(1+\mathcal{V}_{a, b}(s, t)\right)^{-\frac{n+ 2 \overline{r}+2\alpha}{4 \sigma}-1} \|w(s, \cdot)\|_{\overline{r}, 0}, \\
\left\||\xi|^\alpha \mathrm{e}^{-|\xi|^{2 \sigma} \mathcal{V}_{a, b}(s, t)} \widehat{w}(s, \cdot)\right\|_{L^2\left(|\xi| \leq \Lambda_0(t, s)\right)} \lesssim & \left(1+\mathcal{V}_{a, b}(s, t)\right)^{-\frac{n+ 2 \overline{r}+2\alpha}{4 \sigma}} \|w(s, \cdot)\|_{\overline{r}, 0}.
\end{aligned}
$$
From the above estimates, we conclude the proof of \eqref{eq:3.17}.
\end{proof}
For $t>s\geq t_0$, we combine \eqref{eq:3.16} with
the hyperbolic and reduced estimates as in the proof
of \eqref{eq:3.15_L2}.
By \textbf{(B2)}, \textbf{(C2)} and \eqref{eq:exp-parabolic},
for every $c>0$ we can choose $c'\in(0,c)$ such that
$$
\exp\left(-c\int_s^t b(\tau) d\tau\right) \lesssim \frac{\chi(t)}{\chi(s)b(s)} \min\left\{b(t), \frac{a^2(t)}{b(t)\mathcal{V}_{a,b}(s,t)}\right\} \exp\left(-c'\int_s^t b(\tau) d\tau\right).
$$
We can therefore include the exponentially decaying
terms in both estimates.

For $\mathcal{V}_{a,b}(s,t)\geq1$, we apply
Lemma \ref{lemma3.8} on the elliptic part.
For $0<\mathcal{V}_{a,b}(s,t)\leq1$, we use
$$
\sup_{\xi\in\mathbb{R}^n}
|\xi|^{2\sigma}
\exp\left(-|\xi|^{2\sigma}\mathcal{V}_{a,b}(s,t)\right)
\lesssim \mathcal{V}_{a,b}(s,t)^{-1}
$$
and the $L^2$ bound. For frequencies crossing a zone boundary, we compose the corresponding fundamental solutions as above. Using \textbf{(B2)} and \eqref{eq:exp-parabolic}, we obtain two bounds with the same decay factors and the respective coefficients
\begin{equation}\label{eq:ut_coefficients}
b(t) \quad \text{ and } \quad \frac{m^2(t)}{b(t)}+\frac{a^2(t)}{b(t)\mathcal{V}_{a,b}(s,t)}.
\end{equation}
For $t>s$, using $\min\{1,x^{-1}\}\leq2/(1+x)$ for $x>0$, we obtain
\begin{equation}\label{eq:ut_min_bound}
\min\left\{b(t),\frac{m^2(t)}{b(t)}+\frac{a^2(t)}{b(t)\mathcal{V}_{a,b}(s,t)}\right\}\lesssim \frac{m^2(t)}{b(t)}+b(t)\left(1+\frac{b^2(t)}{a^2(t)}\mathcal{V}_{a,b}(s,t)\right)^{-1}.
\end{equation}
The estimates for $s<t_0$ follow from the evolution on bounded time intervals and the composition property.

\begin{proof}[Proof of Theorem \ref{theorem:1.1}]
We write the solution to \eqref{eq:1.1} as
$$
u(t,\cdot)=\varphi_0(t,0,\cdot)*_{(x)}u_0+\varphi_1(t,0,\cdot)*_{(x)}u_1.
$$
For $0\leq\alpha\leq\sigma$, applying \eqref{eq:3.15} with $s=0$ to $u_1$ and the same estimates without $h^{-1}(0,\xi)$ to $u_0$, we obtain
$$
\|u(t,\cdot)\|_{\dot{H}^\alpha}\lesssim \chi(t)\left\{\left(1+\mathcal{V}_{a,b}(0,t)\right)^{-\frac{n+2r_0+2\alpha}{4\sigma}}\|u_0\|_{r_0,\alpha}+\left(1+\mathcal{V}_{a,b}(0,t)\right)^{-\frac{n+2r_1+2\alpha}{4\sigma}}\|u_1\|_{r_1,0}\right\}.
$$
This proves \eqref{estimates_linear} for $j=0$.

For $j=1$ and $\mathcal{V}_{a,b}(0,t)\geq1$, \eqref{eq:3.16}, Lemma \ref{lemma3.8}, the hyperbolic and reduced estimates, together with \eqref{eq:exp-parabolic}, give
$$
\begin{aligned}
\|u_t(t,\cdot)\|_{\dot{H}^\alpha} \lesssim&\ \frac{\chi(t)\left(a^2(t)\left(1+\mathcal{V}_{a,b}(0,t)\right)^{-1}+m^2(t)\right)}{b(t)}\\
&\times\left\{\left(1+\mathcal{V}_{a,b}(0,t)\right)^{-\frac{n+2r_0+2\alpha}{4\sigma}}\|u_0\|_{r_0,\alpha+\sigma}+\left(1+\mathcal{V}_{a,b}(0,t)\right)^{-\frac{n+2r_1+2\alpha}{4\sigma}}\|u_1\|_{r_1,\alpha}\right\}.
\end{aligned}
$$
The case $\mathcal{V}_{a,b}(0,t)\leq1$ follows from the $L^2$ estimates on bounded time intervals. Thus, \eqref{estimates_linear} holds for $j=1$ and $0\leq\alpha\leq\sigma$.

For $j=0$ and $\alpha>\sigma$, we apply the estimate with $\alpha=\sigma$ to $|D|^{\alpha-\sigma}u$. For $j=1$ and $\alpha>\sigma$, we apply the estimate with $\alpha=0$ to $|D|^\alpha u$. In both cases, we use
$$
\Theta_{r_i+\theta}\left(|D|^\theta u_i\right)\leq\Theta_{r_i}(u_i),\quad i=0,1,
$$
with $\theta=\alpha-\sigma$ and $\theta=\alpha$, respectively. The corresponding Sobolev norms are controlled by the data norms in \eqref{estimates_linear}. We therefore obtain \eqref{estimates_linear} for every $\alpha\geq0$ and $j=0,1$. This completes the proof of Theorem \ref{theorem:1.1}.
\end{proof}

\begin{proof}[Proof of Proposition \ref{proposition:5.3}]
Let $u_0(s,\cdot)=0$. Then
$$
u(t,\cdot)=\varphi_1(t,s,\cdot)*_{(x)}u_1(s,\cdot).
$$
For $j=0$, applying \eqref{eq:3.15_L2} with $w=u_1(s,\cdot)$, we obtain
$$
\|u(t,\cdot)\|_{\dot{H}^\alpha} \lesssim \frac{\chi(t)}{\chi(s)b(s)} \left(\left(\frac{A(s)}{\mu(s)}\right)^2+\mathcal{V}_{a,b}(s,t)\right)^{-\frac{\alpha}{2\sigma}} \|u_1(s,\cdot)\|_{L^2},
$$
which proves \eqref{estimates_linear_u_0=0_L2} for $j=0$. If $P_{r_1}(u_1(s,\cdot))_+<\infty$, then \eqref{eq:3.15} gives
$$
\begin{aligned}
\|u(t,\cdot)\|_{\dot{H}^\alpha} \lesssim&\ \frac{\chi(t)}{\chi(s)b(s)}\Bigg\{\left(1+\mathcal{V}_{a,b}(s,t)\right)^{-\frac{n+2r_1}{4\sigma}}\left(\left(\frac{A(s)}{\mu(s)}\right)^2+\mathcal{V}_{a,b}(s,t)\right)^{-\frac{\alpha}{2\sigma}}\\
&+\left(\frac{\mu(s)}{A(s)}\right)^{\frac{\alpha}{\sigma}} \exp\left(-c\int_s^t b(\tau)d\tau\right)\Bigg\} \|u_1(s,\cdot)\|_{r_1,0}.
\end{aligned}
$$
Thus, \eqref{estimates_linear_u_0=0} follows for $j=0$.

For $j=1$, applying the $L^2$-$L^2$ multiplier estimates to $|D|^\alpha u_1(s,\cdot)$ and using \eqref{eq:ut_min_bound}, we obtain \eqref{estimates_linear_u_0=0_L2}.

For $\mathcal{V}_{a,b}(s,t)\geq1$, the elliptic estimate \eqref{eq:3.17}, together with the hyperbolic and reduced estimates and \eqref{eq:ut_min_bound}, gives
$$
\begin{aligned}
\|u_t(t,\cdot)\|_{\dot{H}^\alpha} \lesssim&\ \frac{\chi(t)}{\chi(s)b(s)}\left(\frac{m^2(t)}{b(t)}+b(t)\left(1+\frac{b^2(t)}{a^2(t)}\mathcal{V}_{a,b}(s,t)\right)^{-1}\right)\\
&\times\Bigg\{\left(1+\mathcal{V}_{a,b}(s,t)\right)^{-\frac{n+2r_1+2\alpha}{4\sigma}}+\exp\left(-c\int_s^t b(\tau)d\tau\right)\Bigg\}\|u_1(s,\cdot)\|_{r_1,\alpha}.
\end{aligned}
$$
For $\mathcal{V}_{a,b}(s,t)\leq1$, we have
$$
1+\mathcal{V}_{a,b}(s,t)\sim 1 \quad \text{ and } \quad \|u_1(s,\cdot)\|_{H^\alpha}\lesssim\|u_1(s,\cdot)\|_{r_1,\alpha}.
$$
Hence, \eqref{estimates_linear_u_0=0_L2} yields \eqref{estimates_linear_u_0=0} for $j=1$. This completes the proof of Proposition \ref{proposition:5.3}.
\end{proof}

\section{Existence of global (in time) solutions} \label{section_5}
\subsection{Outline of the proof}

Let $\varphi_i(t,s,x), i=0,1$, denote the fundamental solutions associated with the linear problem \eqref{eq:3.1}. Then 
$$
u^{\operatorname{lin}}(t,x):= \varphi_0(t, 0, x) *_{(x)} u_0(x)+\varphi_1(t, 0, x) *_{(x)} u_1(x)
$$ 
solves the corresponding linear problem \eqref{eq:1.1}. A function $u\in\mathcal{C}([0,T],H^\alpha)$ is called a mild solution to \eqref{eq:1.0} on $[0,T]$ if
$$
u(t,x) = \varphi_0(t,0,x)*_{(x)}u_0(x) + \varphi_1(t,0,x)*_{(x)}u_1(x) + \int_0^t \varphi_1(t,s,x)*_{(x)}|u(s,x)|^p d s := u^{\operatorname{lin}}(t, x)+u^{\operatorname{nl}}(t, x)
$$
for any $t\in[0,T]$. For any $T>0$, we define the operator $\mathcal{N}: u \in \mathcal{X}(T) \mapsto \mathcal{N} u=\mathcal{N} u(t, x) := u^{\operatorname{lin}}(t, x)+u^{\operatorname{nl}}(t, x)$, where $\mathcal{X}(T)$ is the evolution space introduced later. It remains to show that $\mathcal{N}$ is a contraction on a suitable ball of $\mathcal{X}(T)$. For this purpose, we are going to prove that
\begin{align} \label{ineq:4.1}
\|\mathcal{N} u\|_{\mathcal{X}(T)} \lesssim & \|(u_0, u_1)\|_{\mathcal{D}_{r_0, r_1}^{\alpha,0}}+\|u\|_{\mathcal{X}(T)}^p, \\
\label{ineq:4.2}
\|\mathcal{N} u-\mathcal{N} v\|_{\mathcal{X}(T)} \lesssim & \|u-v\|_{\mathcal{X}(T)}\left(\|u\|_{\mathcal{X}(T)}^{p-1}+\|v\|_{\mathcal{X}(T)}^{p-1}\right).
\end{align}
If $\|(u_0, u_1)\|_{\mathcal{D}_{r_0, r_1}^{\alpha,0}}$ is sufficiently small, then \eqref{ineq:4.1} and \eqref{ineq:4.2} give the existence of a unique global (in time) mild solution by Banach's fixed point theorem.
\begin{lemma}\label{lemma:split_linear}
Assume \textbf{(A)}--\textbf{(C)}. Let $1\leq\omega\leq2$, $0\leq\ell\leq\sigma$, and put
$$
d=\frac{n}{2\sigma}\left(\frac1\omega-\frac12\right).
$$
For $g\in L^\omega\cap L^2$ and $\mathcal{V}_{a,b}(s,t)\geq1$, we have
\begin{align}\label{eq:split_linear}
\| |D|^\ell\varphi_1(t,s,\cdot)*_{(x)}g\|_{L^2} \lesssim&\frac{\chi(t)}{\chi(s)b(s)}\bigg\{ \left(1+\mathcal{V}_{a,b}(s,t)\right)^{-d} \left(\left(\frac{A(s)}{\mu(s)}\right)^2+\mathcal{V}_{a,b}(s,t)\right)^{-\frac{\ell}{2\sigma}} \|g\|_{L^\omega\cap L^2}\notag\\
&+\left(\frac{\mu(s)}{A(s)}\right)^{\frac{\ell}{\sigma}} \exp\left(-c\int_s^t b(\tau)d\tau\right)\|g\|_{L^2} \bigg\}.
\end{align}
\end{lemma}
\begin{proof}
Let $q\in[2,\infty]$ satisfy $1/q=1/\omega-1/2$. In the elliptic zone, we have $|\xi|^\sigma\lesssim\mu(s)/A(s)$ and
$$
|\xi|^\ell\exp\left(-c|\xi|^{2\sigma}\mathcal{V}_{a,b}(s,t)\right)\lesssim\left(\left(\frac{A(s)}{\mu(s)}\right)^2+\mathcal{V}_{a,b}(s,t)\right)^{-\frac{\ell}{2\sigma}}\exp\left(-c'|\xi|^{2\sigma}\mathcal{V}_{a,b}(s,t)\right)
$$
for some $0<c'<c$. By H\"older's and Hausdorff-Young's inequalities, we obtain
$$
\left\||\xi|^\ell\widehat{\varphi}_1(t,s,\xi)\widehat{g}(\xi)\right\|_{L^2(|\xi|\leq\Lambda_0(t,s))}\lesssim\frac{\chi(t)}{\chi(s)b(s)}\left(1+\mathcal{V}_{a,b}(s,t)\right)^{-d}\left(\left(\frac{A(s)}{\mu(s)}\right)^2+\mathcal{V}_{a,b}(s,t)\right)^{-\frac{\ell}{2\sigma}}\|g\|_{L^\omega}.
$$
For $|\xi|\geq\Lambda_\infty(t,s)$, the estimates in the proof of \eqref{eq:3.15} give
\begin{equation}\label{eq:split-high}
|\xi|^\ell|\widehat{\varphi}_1(t,s,\xi)|\lesssim\frac{\chi(t)}{\chi(s)b(s)}\left(\frac{\mu(s)}{A(s)}\right)^{\frac{\ell}{\sigma}}\exp\left(-c\int_s^t b(\tau)d\tau\right).
\end{equation}
For $\Lambda_0(t,s)<|\xi|<\Lambda_\infty(t,s)$, we take $t_\xi\in[s,t]$ with $\Lambda(t_\xi)=|\xi|$. If $\mu(t)/A(t)$ is non-decreasing, we obtain
$$
\left|\widehat{\varphi}_1(t,s,\xi)\right| \lesssim \frac{\chi(t)}{\chi(s)b(s)} \exp\left(-c\int_s^{t_\xi}b(\tau)d\tau-|\xi|^{2\sigma}\mathcal{V}_{a,b}(t_\xi,t)\right).
$$
By \textbf{(B2)}, there exists $C>0$ such that
$$
\frac{\mu(t_\xi)}{A(t_\xi)} \frac{A(s)}{\mu(s)} \leq \left(\frac{A(t_\xi)}{A(s)}\right)^C.
$$
Since $\mu(t)\to\infty$, we have
\begin{equation}\label{eq:intermediate-absorption}
\left(\frac{A(t_\xi)}{A(s)}\right)^L\exp\left(-\frac{c}{4}\int_s^{t_\xi}b(\tau)d\tau\right)=\exp\left(\int_s^{t_\xi}\left(L-\frac{c}{4}\mu(\tau)\right)\frac{a(\tau)}{A(\tau)}d\tau\right)\lesssim 1 \quad \text{ for every } L>0.
\end{equation}
From $|\xi|^\sigma\sim\mu(t_\xi)/A(t_\xi)$, we obtain
$$
1+|\xi|^{2\sigma}\mathcal{V}_{a,b}(s,t_\xi) \lesssim \left(\frac{A(t_\xi)}{A(s)}\right)^{2C} \left( 1+\left(\frac{\mu(s)}{A(s)}\right)^2 \mathcal{V}_{a,b}(s,t_\xi) \right).
$$
By \eqref{eq:exp-parabolic} and \eqref{eq:intermediate-absorption} with $L=2CN$, we have
$$
\begin{aligned}
\exp\left(-\frac{c}{2}\int_s^{t_\xi}b(\tau)d\tau\right) \lesssim & C_N\left(\frac{A(t_\xi)}{A(s)}\right)^{-2CN}\left(1+\left(\frac{\mu(s)}{A(s)}\right)^2\mathcal{V}_{a,b}(s,t_\xi)\right)^{-N}\\
\lesssim & C_N \left(1+|\xi|^{2\sigma}\mathcal{V}_{a,b}(s,t_\xi)\right)^{-N},
\end{aligned}
$$
which gives
$$
\begin{aligned}
\exp\left(-\frac{c}{2}\int_s^{t_\xi}b(\tau)d\tau-\frac{1}{2}|\xi|^{2\sigma}\mathcal{V}_{a,b}(t_\xi,t)\right) \lesssim &C_N \left(1+|\xi|^{2\sigma}\mathcal{V}_{a,b}(s,t_\xi)\right)^{-N}\left(1+|\xi|^{2\sigma}\mathcal{V}_{a,b}(t_\xi,t)\right)^{-N}\\
\lesssim &C_N \left(1+|\xi|^{2\sigma}\mathcal{V}_{a,b}(s,t)\right)^{-N}.
\end{aligned}
$$
If $\mu(t)/A(t)$ is non-increasing, one gets
$$
\left|\widehat{\varphi}_1(t,s,\xi)\right| \lesssim \frac{\chi(t)}{\chi(s)b(s)}\exp\left(-|\xi|^{2\sigma}\mathcal{V}_{a,b}(s,t_\xi)-c\int_{t_\xi}^t b(\tau)d\tau\right).
$$
By \eqref{eq:exp-parabolic} and $|\xi|^\sigma\sim\mu(t_\xi)/A(t_\xi)$, we obtain
$$
\begin{aligned}
\exp\left(-\frac{1}{2}|\xi|^{2\sigma}\mathcal{V}_{a,b}(s,t_\xi)-\frac{c}{2}\int_{t_\xi}^t b(\tau)d\tau\right) \lesssim &C_N \left(1+|\xi|^{2\sigma}\mathcal{V}_{a,b}(s,t_\xi)\right)^{-N}\left(1+|\xi|^{2\sigma}\mathcal{V}_{a,b}(t_\xi,t)\right)^{-N}\\
\lesssim &C_N\left(1+|\xi|^{2\sigma}\mathcal{V}_{a,b}(s,t)\right)^{-N}.
\end{aligned}
$$
In both cases, splitting the exponential factors and using the estimates in the proof of \eqref{eq:3.15_L2}, we obtain
$$
|\xi|^\ell\left|\widehat{\varphi}_1(t,s,\xi)\right|\lesssim C_N\frac{\chi(t)}{\chi(s)b(s)}\left(\left(\frac{A(s)}{\mu(s)}\right)^2+\mathcal{V}_{a,b}(s,t)\right)^{-\frac{\ell}{2\sigma}}\left(1+|\xi|^{2\sigma}\mathcal{V}_{a,b}(s,t)\right)^{-N} \quad \text{ for every } N>0.
$$
Choose $N>d$. By scaling, we derive
$$
\left\| \left( 1+|\xi|^{2\sigma}\mathcal{V}_{a,b}(s,t) \right)^{-N} \right\|_{L^q} \lesssim \mathcal{V}_{a,b}(s,t)^{-d}.
$$
Since $\mathcal{V}_{a,b}(s,t)\geq1$, by H\"older's and Hausdorff-Young's inequalities, one finds
\begin{equation}\label{eq:split-intermediate}
\left\||\xi|^\ell\widehat{\varphi}_1(t,s,\xi)\widehat{g}(\xi)\right\|_{L^2(\Lambda_0(t,s)<|\xi|<\Lambda_\infty(t,s))}\lesssim \frac{\chi(t)}{\chi(s)b(s)}\left(1+\mathcal{V}_{a,b}(s,t)\right)^{-d}\left(\left(\frac{A(s)}{\mu(s)}\right)^2+\mathcal{V}_{a,b}(s,t)\right)^{-\frac{\ell}{2\sigma}}\|g\|_{L^\omega}.
\end{equation}
Therefore, combining the elliptic estimate with \eqref{eq:split-high} and \eqref{eq:split-intermediate} we may conclude \eqref{eq:split_linear}.
\end{proof}

\subsection{Proof of Theorem \ref{theorem:1.2}}
We put $\mathcal{A}_1(p):=n/(2\omega\sigma)(p-1)$. By Proposition \ref{proposition:inf_beta}, there exists $\delta>0$ such that
\begin{equation}\label{eq:nonlinear_weight_chi}
\frac{\chi(t)^{p-1}}{a^2(t)} \left(1+\mathcal{V}_{a,b}(0,t)\right)^{-\mathcal{A}_1(p)} \lesssim \left(1+\mathcal{V}_{a,b}(0,t)\right)^{-1-\delta}, \quad t\geq0.
\end{equation}
For $T>0$, we define $\mathcal{X}(T):=\mathcal{C}([0,T],H^\alpha)$, endowed with the norm
$$
\begin{aligned}
\|u\|_{\mathcal{X}(T)} := \sup_{t\in[0,T]} \chi(t)^{-1} \Big(&\left(1+\mathcal{V}_{a,b}(0,t)\right)^{\frac{n+2\min\{r_0,r_1\}}{4\sigma}}\|u(t,\cdot)\|_{L^2} \\
&+\left(1+\mathcal{V}_{a,b}(0,t)\right)^{ \frac{n+2\min\{r_0,r_1\}+2\alpha}{4\sigma}} \big\||D|^\alpha u(t,\cdot)\big\|_{L^2} \Big).
\end{aligned}
$$
\medskip
Since $\min\{r_0,r_1\}\in[-n/2,0]$, we have $1\leq\omega\leq2$. Applying the fractional Gagliardo-Nirenberg inequality from Lemma \ref{Gagliardo} with $\theta_1:= n(p-1)/(2 \alpha p)$ and $\theta_2:= n(\omega p-2)/(2 \omega \alpha p)$ together with the definition of the space $\mathcal{X}(T)$, we derive
\begin{align}\label{ineq:4.3}
\||u(\tau,\cdot)|^p\|_{L^2} &=\|u(\tau,\cdot)\|_{L^{2p}}^p \lesssim\chi(\tau)^p \left(1+\mathcal{V}_{a,b}(0,\tau)\right)^{ -\frac{n+\min\{r_0,r_1\}}{2\sigma}p +\frac{n}{4\sigma}} \|u\|_{\mathcal{X}(T)}^p,\\
\||u(\tau,\cdot)|^p\|_{L^\omega} &=\|u(\tau,\cdot)\|_{L^{\omega p}}^p \lesssim\chi(\tau)^p \left(1+\mathcal{V}_{a,b}(0,\tau)\right)^{ -\frac{n+\min\{r_0,r_1\}}{2\sigma}p +\frac{n}{2\omega\sigma}} \|u\|_{\mathcal{X}(T)}^p.
\label{ineq:4.4}
\end{align}
provided that $p \in [2/\omega, \infty)$ if $n \leq 2 \alpha$, or $p \in [2/\omega, n/(n-2 \alpha)]$ if $n>2\alpha$ and $n(2-\omega)\leq4\alpha$. We split the Duhamel integral into two parts. Since $\mathcal{V}_{a,b}(0,t)$ is strictly increasing in $t$ and $\mathcal{V}_{a,b}(0,0)=0$, for every $t>0$ there is a unique $s_t\in(0,t)$ such that
$$
\mathcal{V}_{a,b}(0,s_t)=\frac12\mathcal{V}_{a,b}(0,t).
$$
Thus, it follows that
\begin{equation}\label{ineq:V_a,b}
\mathcal{V}_{a,b}(s,t)\geq\frac12\mathcal{V}_{a,b}(0,t)\quad\text{ if } 0\leq s\leq s_t, \quad \left(1+\mathcal{V}_{a,b}(0,s)\right)\geq\frac12\left(1+\mathcal{V}_{a,b}(0,t)\right)\quad \text{ if }s_t\leq s\leq t.
\end{equation}
For $\mathcal{V}_{a,b}(0,t)\leq2$, the $L^2$-$L^2$ estimate in Corollary \ref{corollary:linear} and \eqref{ineq:4.3} give
$$
\||D|^\ell u^{\operatorname{nl}}(t,\cdot)\|_{L^2} \lesssim\|u\|_{\mathcal{X}(T)}^p, \quad \ell=0,\alpha.
$$
By \textbf{(B4)}, these times are bounded, and we have $\chi(t)\gtrsim1$ and $1+\mathcal{V}_{a,b}(0,t)\lesssim1$. Hence, we get
$$
\||D|^\ell u^{\operatorname{nl}}(t,\cdot)\|_{L^2} \lesssim\chi(t) \left(1+\mathcal{V}_{a,b}(0,t)\right)^{ -\frac{n+2\min\{r_0,r_1\}+2\ell}{4\sigma}} \|u\|_{\mathcal{X}(T)}^p.
$$
For $\mathcal{V}_{a,b}(0,t)\geq2$, by \eqref{ineq:V_a,b}, we have $\mathcal{V}_{a,b}(s,t)\geq1$ for $s\in[0,s_t]$. Using Lemma \ref{lemma:split_linear} on $[0,s_t]$, Corollary \ref{corollary:linear} on $[s_t,t]$ and \eqref{ineq:4.3}-\eqref{ineq:4.4}, we obtain
$$
\||D|^\ell u^{\operatorname{nl}}(t)\|_{L^2} \lesssim\|u\|_{\mathcal{X}(T)}^p (\mathcal{H}_\ell+\mathcal{E}_\ell+\mathcal{R}_\ell), \quad \ell=0,\alpha,
$$
where
$$
\begin{aligned}
\mathcal{H}_\ell:=&\chi(t)\int_0^{s_t} \frac{ \chi(s)^{p-1} \left(1+\mathcal{V}_{a,b}(0,s)\right)^{-\mathcal{A}_1(p)} }{b(s)} (1+\mathcal{V}_{a,b}(s,t))^{ -\frac{n+2\min\{r_0,r_1\}}{4\sigma}} \left( \left(\frac{A(s)}{\mu(s)}\right)^2 +\mathcal{V}_{a,b}(s,t) \right)^{-\frac{\ell}{2\sigma}}ds,\\
\mathcal{E}_\ell:=&\chi(t)\int_0^{s_t} \frac{ \chi(s)^{p-1} \left(1+\mathcal{V}_{a,b}(0,s)\right)^{ -\mathcal{A}_1(p)-\frac{n+2\min\{r_0,r_1\}}{4\sigma}} }{b(s)} \left(\frac{\mu(s)}{A(s)}\right)^{\frac{\ell}{\sigma}} \exp\left(-c\int_s^t b(\tau)d\tau\right)ds,\\
\mathcal{R}_\ell:=&\chi(t)\int_{s_t}^t \frac{ \chi(s)^{p-1} \left(1+\mathcal{V}_{a,b}(0,s)\right)^{-\mathcal{A}_2(p)} }{b(s)} \left( \left(\frac{A(s)}{\mu(s)}\right)^2 +\mathcal{V}_{a,b}(s,t) \right)^{-\frac{\ell}{2\sigma}}ds,
\end{aligned}
$$
with
$$
\mathcal{A}_2(p):= \frac{n}{2\omega\sigma}p-\frac{n}{4\sigma} = \mathcal{A}_1(p)+\frac{n+2\min\{r_0,r_1\}}{4\sigma}.
$$
\begin{itemize}[leftmargin=*]
\item[$\bullet$] \textbf{Estimate of $\mathcal{H}_\ell$:} For $s\leq s_t$ and $\mathcal{V}_{a,b}(0,t)\geq2$, we have $1+\mathcal{V}_{a,b}(s,t)\gtrsim 1+\mathcal{V}_{a,b}(0,t)$ and $\mathcal{V}_{a,b}(s,t)\gtrsim 1+\mathcal{V}_{a,b}(0,t)$. Hence,
$$
\mathcal{H}_\ell\lesssim\chi(t) \left(1+\mathcal{V}_{a,b}(0,t)\right)^{ -\frac{n+2\min\{r_0,r_1\}+2\ell}{4\sigma}} \int_0^{s_t} \frac{ \chi(s)^{p-1} \left(1+\mathcal{V}_{a,b}(0,s)\right)^{-\mathcal{A}_1(p)} }{b(s)}ds.
$$
By \eqref{eq:nonlinear_weight_chi}, using the change of variables $r=\mathcal{V}_{a,b}(0,s)$ we derive
$$
\begin{aligned}
\int_0^{s_t} \frac{ \chi(s)^{p-1} \left(1+\mathcal{V}_{a,b}(0,s)\right)^{-\mathcal{A}_1(p)} }{b(s)}ds= &\int_0^{s_t} \frac{ \chi(s)^{p-1} \left(1+\mathcal{V}_{a,b}(0,s)\right)^{-\mathcal{A}_1(p)} }{a^2(s)} d\mathcal{V}_{a,b}(0,s) \\
\lesssim &\int_0^{\mathcal{V}_{a,b}(0,t)/2} (1+r)^{-1-\delta}dr \lesssim1,
\end{aligned}
$$
which entails
\begin{equation}\label{eq:nonlinear_H_bound}
\mathcal{H}_\ell\lesssim\chi(t) \left(1+\mathcal{V}_{a,b}(0,t)\right)^{ -\frac{n+2\min\{r_0,r_1\}+2\ell}{4\sigma}}.
\end{equation}

\item[$\bullet$] \textbf{Estimate of $\mathcal{E}_\ell$:} From \textbf{(B2)}, we have
$$
\left(\frac{A^2(s)}{\mu(s)}\right)' \leq(2+M)\frac{a^2(s)}{b(s)},
$$
that is, $A^2(s)/\mu(s)\lesssim 1+\mathcal{V}_{a,b}(0,s)$. Together with $\inf_{s\geq0}\mu(s)>0$, we obtain
$$
\left(\frac{\mu(s)}{A(s)}\right)^2 \left(1+\mathcal{V}_{a,b}(0,s)\right) \gtrsim\mu(s)\gtrsim1.
$$
We put $N=(n+2\min\{r_0,r_1\}+2\ell)/(4\sigma)$. For $N>0$, applying \eqref{eq:exp-parabolic}, we obtain
$$
\begin{aligned}
&\left(\frac{\mu(s)}{A(s)}\right)^{\frac{\ell}{\sigma}} \exp\left(-c\int_s^t b(\tau)d\tau\right) \left(1+\mathcal{V}_{a,b}(0,s)\right)^{ -\frac{n+2\min\{r_0,r_1\}}{4\sigma}}\\
&\quad \lesssim\left(\frac{\mu(s)}{A(s)}\right)^{\frac{\ell}{\sigma}} \left( \left(\frac{\mu(s)}{A(s)}\right)^2 \mathcal{V}_{a,b}(s,t) \right)^{-N} \left(1+\mathcal{V}_{a,b}(0,s)\right)^{ -\frac{n+2\min\{r_0,r_1\}}{4\sigma}}\\
&\quad =\left( \left(\frac{\mu(s)}{A(s)}\right)^2 \left(1+\mathcal{V}_{a,b}(0,s)\right) \right)^{-\frac{n+2\min\{r_0,r_1\}}{4\sigma}} \mathcal{V}_{a,b}(s,t)^{-N}\lesssim \left(1+\mathcal{V}_{a,b}(0,t)\right)^{ -\frac{n+2\min\{r_0,r_1\}}{4\sigma} -\frac{\ell}{2\sigma}}
\end{aligned}
$$
for $s\in[0,s_t]$. If $N=0$, then the same estimate follows from $\exp\left(-c\int_s^t b(\tau)d\tau\right)\leq1$. Hence, we obtain
\begin{align}\label{eq:nonlinear_E_bound}
\mathcal{E}_\ell \lesssim & \chi(t) \left(1+\mathcal{V}_{a,b}(0,t)\right)^{ -\frac{n+2\min\{r_0,r_1\}+2\ell}{4\sigma}} \int_0^{s_t} \frac{ \chi(s)^{p-1} \left(1+\mathcal{V}_{a,b}(0,s)\right)^{-\mathcal{A}_1(p)} }{b(s)}ds \notag\\
\lesssim & \chi(t) \left(1+\mathcal{V}_{a,b}(0,t)\right)^{ -\frac{n+2\min\{r_0,r_1\}+2\ell}{4\sigma}}.
\end{align}

\item[$\bullet$] \textbf{Estimate of $\mathcal{R}_\ell$:} For $s_t\leq s\leq t$, we have $1+\mathcal{V}_{a,b}(0,s)\sim 1+\mathcal{V}_{a,b}(0,t)$. From $\ell/(2\sigma)\in[0,1/2]$ and \eqref{eq:nonlinear_weight_chi}, it yields
$$
\begin{aligned}
\mathcal{R}_\ell \lesssim & \chi(t)\int_{s_t}^t \frac{\chi(s)^{p-1}}{a^2(s)} \left(1+\mathcal{V}_{a,b}(0,s)\right)^{-\mathcal{A}_2(p)} \mathcal{V}_{a,b}(s,t)^{-\frac{\ell}{2\sigma}} d\mathcal{V}_{a,b}(0,s)\\
\lesssim & \chi(t) \left(1+\mathcal{V}_{a,b}(0,t)\right)^{ -\frac{n+2\min\{r_0,r_1\}}{4\sigma}-1-\delta} \int_0^{\mathcal{V}_{a,b}(0,t)/2} r^{-\frac{\ell}{2\sigma}}dr\\
\lesssim & \chi(t) \left(1+\mathcal{V}_{a,b}(0,t)\right)^{ -\frac{n+2\min\{r_0,r_1\}+2\ell}{4\sigma}-\delta}\lesssim\chi(t) \left(1+\mathcal{V}_{a,b}(0,t)\right)^{ -\frac{n+2\min\{r_0,r_1\}+2\ell}{4\sigma}},
\end{aligned}
$$
\end{itemize}

Combining the estimates for $\mathcal{H}_\ell$, $\mathcal{E}_\ell$ and $\mathcal{R}_\ell$, we obtain
$$
\||D|^\ell u^{\operatorname{nl}}(t)\|_{L^2} \lesssim\chi(t) \left(1+\mathcal{V}_{a,b}(0,t)\right)^{-\frac{n+2\min\{r_0,r_1\}+2\ell}{4\sigma}}\|u\|_{\mathcal{X}(T)}^p, \quad \ell=0,\alpha.
$$
Hence, we get
\begin{equation}\label{eq:nonlinear_X_bound}
\|u^{\operatorname{nl}}\|_{\mathcal{X}(T)} \lesssim\|u\|_{\mathcal{X}(T)}^p.
\end{equation}

For $u,v\in\mathcal{X}(T)$, it is obvious to see that
$$
u^{\operatorname{nl}}(t)-v^{\operatorname{nl}}(t)=\int_0^t\varphi_1(t,s,\cdot)*_{(x)}\left(|u(s,\cdot)|^p-|v(s,\cdot)|^p\right)ds.
$$
Since $p>1$, we have
$$
\left|| u(s, x)|^p-| v(s, x)|^p \right| \lesssim | u(s, x)- v(s, x)| \left(| u(s, x)|^{p-1}+| v(s, x)|^{p-1}\right).
$$
Hence, H\"{o}lder's inequality implies
$$
\begin{aligned}
\left\| \left| u(s, \cdot)\right|^p-\left| v(s, \cdot)\right|^p\right\|_{L^{\omega}}  \lesssim & \left\| u(s, \cdot)- v(s, \cdot)\right\|_{L^{\omega p}}\left( \left\| u(s, \cdot)\right\|_{L^{\omega p}}^{p-1}+\left\| v(s, \cdot)\right\|_{L^{\omega p}}^{p-1}\right), \\
\left\| \left| u(s, \cdot)\right|^p-\left| v(s, \cdot)\right|^p\right\|_{L^2}  \lesssim & \left\| u(s, \cdot)- v(s, \cdot)\right\|_{L^{2 p}}\left( \left\| u(s, \cdot)\right\|_{L^{2 p}}^{p-1}+\left\| v(s, \cdot)\right\|_{L^{2 p}}^{p-1}\right).
\end{aligned}
$$
Applying Lemma \ref{Gagliardo} to $u-v$, $u$ and $v$, with $q=\omega p$ and $q=2p$, we get
\begin{align}
&\||u(s)|^p-|v(s)|^p\|_{L^\omega} \lesssim \chi(s)^p \left(1+\mathcal{V}_{a,b}(0,s)\right)^{-\frac{n}{2\omega\sigma}(p-1)} \|u-v\|_{\mathcal{X}(T)} \left(\|u\|_{\mathcal{X}(T)}^{p-1} +\|v\|_{\mathcal{X}(T)}^{p-1}\right), \label{eq:difference_Lomega} \\
&\||u(s)|^p-|v(s)|^p\|_{L^2} \lesssim \chi(s)^p \left(1+\mathcal{V}_{a,b}(0,s)\right)^{-\frac{n}{2\omega\sigma}p+\frac{n}{4\sigma}}\|u-v\|_{\mathcal{X}(T)}\left(\|u\|_{\mathcal{X}(T)}^{p-1}+\|v\|_{\mathcal{X}(T)}^{p-1}\right). \label{eq:difference_L2}
\end{align}
For $\mathcal{V}_{a,b}(0,t)\geq2$, we apply the previous estimates for $\mathcal{H}_\ell$, $\mathcal{E}_\ell$ and $\mathcal{R}_\ell$ to \eqref{eq:difference_Lomega} and \eqref{eq:difference_L2}. For $\mathcal{V}_{a,b}(0,t)\leq2$, we use the $L^2$-$L^2$ estimate in Corollary \ref{corollary:linear}. Hence,
\begin{equation}\label{eq:nonlinear_difference_X_bound}
\|u^{\operatorname{nl}}-v^{\operatorname{nl}}\|_{\mathcal{X}(T)} \lesssim \|u-v\|_{\mathcal{X}(T)}\left(\|u\|_{\mathcal{X}(T)}^{p-1}+\|v\|_{\mathcal{X}(T)}^{p-1}\right).
\end{equation}
It remains to verify that $u^{\operatorname{nl}}\in\mathcal{C}([0,T],H^\alpha)$. By \eqref{ineq:4.3} and the definition of $\mathcal{X}(T)$, it holds $\sup_{s\in[0,T]} \left\||u(s,\cdot)|^p\right\|_{L^2} <\infty$. Then, we apply the $L^2$-$L^2$ estimate in Corollary \ref{corollary:linear} with $j=0$, for $\ell=0$ and $\ell=\alpha$, to obtain
$$
\left\| \varphi_1(t,s,\cdot)*_{(x)}|u(s,\cdot)|^p \right\|_{H^\alpha} \leq C_T \left\||u(s,\cdot)|^p\right\|_{L^2}, \quad 0\leq s\leq t\leq T,
$$
where $C_T>0$ is independent of $s$ and $t$. Fix $t_0\in[0,T]$. For $t\geq t_0$, we write
$$
u^{\operatorname{nl}}(t)-u^{\operatorname{nl}}(t_0)=\int_0^{t_0} \left(\varphi_1(t,s,\cdot)-\varphi_1(t_0,s,\cdot) \right) *_{(x)}|u(s,\cdot)|^p ds+\int_{t_0}^{t}\varphi_1(t,s,\cdot) *_{(x)}|u(s,\cdot)|^p ds.
$$
By the strong continuity of the linear evolution operator and the dominated convergence theorem, we get
$$
\int_0^{t_0}\left\|\left(\varphi_1(t,s,\cdot)-\varphi_1(t_0,s,\cdot)\right) *_{(x)}|u(s,\cdot)|^p\right\|_{H^\alpha} ds \longrightarrow 0 \quad \text{ as } t\to t_0,
$$
which gives
$$
\left\|\int_{t_0}^{t} \varphi_1(t,s,\cdot)*_{(x)}|u(s,\cdot)|^p ds \right\|_{H^\alpha} \lesssim \int_{t_0}^{t} \left\||u(s,\cdot)|^p\right\|_{L^2}ds \longrightarrow 0 \quad \text{ as } t\to t_0.
$$
The case $t<t_0$ follows in the same way. Hence, $u^{\operatorname{nl}}\in\mathcal{C}([0,T],H^\alpha)$.

By Theorem \ref{theorem:1.1}, we get
$$
\begin{aligned}
\|u^{\operatorname{lin}}(t,\cdot)\|_{L^2} \lesssim&\ \chi(t) \left(1+\mathcal{V}_{a,b}(0,t)\right)^{ -\frac{n+2\min\{r_0,r_1\}}{4\sigma}} \left\|(u_0,u_1)\right\|_{ \mathcal{D}_{r_0,r_1}^{\alpha,0}},\\
\||D|^\alpha u^{\operatorname{lin}}(t,\cdot)\|_{L^2} \lesssim&\ \chi(t) \left(1+\mathcal{V}_{a,b}(0,t)\right)^{ -\frac{n+2\min\{r_0,r_1\}+2\alpha}{4\sigma}} \left\|(u_0,u_1)\right\|_{ \mathcal{D}_{r_0,r_1}^{\alpha,0}},
\end{aligned}
$$
to conclude that
\begin{equation}\label{eq:linear_X_bound}
\|u^{\operatorname{lin}}\|_{\mathcal{X}(T)} \lesssim \left\|(u_0,u_1)\right\|_{\mathcal{D}_{r_0,r_1}^{\alpha,0}}.
\end{equation}
Combining \eqref{eq:linear_X_bound} and \eqref{eq:nonlinear_X_bound}, we obtain \eqref{ineq:4.1}. From \eqref{eq:nonlinear_difference_X_bound}, we obtain \eqref{ineq:4.2}. This completes the proof of Theorem \ref{theorem:1.2}.


\section*{Acknowledgments}
Duc An Phan sincerely acknowledges the financial support provided by the Banking Academy of Vietnam. The research of Tuan Anh Dao is funded by the Vietnam National Foundation for Science and Technology Development (NAFOSTED) under grant number 101.02-2025.23.
\medskip

\textbf{Conflict of interest:} The authors declare no conflict of interest.\medskip

\textbf{Data Availability Statement:} Our manuscript has no associated data, and no new data were generated during the study.

\appendix

\section{Decay indicators and decay characters} \label{section_Decay_characters}

\begin{definition}[see \cite{Brandolese2016}] \label{definition_2.1}
Let $f \in L^2$ and $r \in [-n/2, \infty)$. The lower and upper decay indicators of $f$ are defined by
$$
\begin{aligned}
P_r(f)_{-}=\liminf _{\rho \rightarrow 0^{+}} \rho^{-2 r-n} \int_{B(\rho)}|\widehat{f}(\xi)|^2 d \xi \quad \text{ and } \quad P_r(f)_{+}=\limsup _{\rho \rightarrow 0^{+}} \rho^{-2 r-n} \int_{B(\rho)}|\widehat{f}(\xi)|^2 d \xi.
\end{aligned}
$$
\end{definition}

\begin{definition}[see \cite{Brandolese2016}]
\label{definition_2.2}
For $f \in L^2$, define the lower and upper decay characters of $f$ by
$$
r(f)_{-}=\inf \left\{r \in \left[-n/2,\infty\right): P_r(f)_{-}>0\right\} \quad \text{ and } \quad r(f)_{+}=\sup \left\{r \in \left[-n/2,\infty\right): P_r(f)_{+}<\infty\right\}.
$$
We use the convention $\inf\varnothing=+\infty$. Then, $r^*=r^*\left(f\right)$ is called the decay character of $f$ if there exists $r^* \in (-n/2, \infty)$ such that
$$
r^*=\min \left\{r \in \left[-n/2,\infty\right): P_r(f)_{-}>0\right\}=\max \left\{r \in \left[-n/2,\infty\right): P_r(f)_{+}<\infty\right\}.
$$
\end{definition}

\begin{remark}
\rm{
For $f \in L^2$, the quantities $r(f)_{-}$ and $r(f)_{+}$ are well-defined and satisfy $-n/2 \leq r(f)_{+} \leq r(f)_{-} \leq\infty$. 
}
\end{remark}

\begin{example} \label{Example:2.0}
\rm{
If $f \in L^\eta \cap L^2$ with $\eta \in [1, 2]$, the Hausdorff-Young and H\"{o}lder inequalities (see \cite{Niche2022}) give
$$
\Theta_r(f)\lesssim\|f\|_{L^\eta}^2 \quad \text{ and hence } \quad P_r(f)_{+}<\infty \quad \text{ for } \quad r=-n\left(1-1/\eta\right).
$$
Consequently, $r(f)_+\geq-n(1-1/\eta)$. If the decay character $r^*(f)$ exists, then $r^*(f)\geq-n(1-1/\eta)$. In particular, if $f \in L^1 \cap L^2$ and there exist $\varepsilon>0$ and constants $C_1, C_2>0$ such that $C_1 \leq|\widehat{f}(\xi)| \leq C_2$ for $|\xi| \leq \varepsilon$, then $r^*\left(f\right)=0$.
}
\end{example}

\begin{example} \label{Example:2.1}
\rm{
Let $L^{1, \eta}:= L^1\left(\mathbb{R}^n,(1+|x|)^\eta d x\right)$ for $\eta \in[0,1]$. If $f \in L^2\cap L^{1, \eta}$ has nonzero mean, then $\widehat{f}$ is continuous and $\widehat{f}(0)\neq0$, and consequently $r^*(f)=0$. If $f$ has zero mean, then
$$
|\widehat{f}(\xi)|=\left|\int_{\mathbb{R}^n}(\mathrm{e}^{-ix\cdot\xi}-1)f(x)dx\right|\lesssim |\xi|^\eta\|f\|_{L^{1,\eta}}.
$$
Thus $\Theta_\eta(f)\lesssim\|f\|_{L^{1,\eta}}^2$ and $P_\eta(f)_+<\infty$, and therefore $r(f)_+\geq\eta$; in particular, $r^*(f)\geq\eta$ whenever the decay character exists.
} 
\end{example}

\begin{example} \label{Example:2.2}
\rm{
We define
$$
\dot{H}^{-\eta}:=\left\{f \in S': \|f\|_{\dot{H}^{-\eta}}:=\big\|(-\Delta)^{-\frac{\eta}{2}} f\big\|_{L^2}=\big\|\mathcal{F}^{-1}\big(|\xi|^{-\eta} \widehat{f}\big)\big\|_{L^2}<\infty\right\},
$$
where $S'\left(\mathbb{R}^n\right)$ is the space of tempered distributions. If $f \in L^2 \cap \dot{H}^{-\eta}$ with $\eta \geq 0$, then
$$
\rho^{-2 r-n} \int_{B(\rho)}\left|\widehat{f}(\xi)\right|^2 d \xi \leq \rho^{2 \eta-2 r-n} \int_{B(\rho)}|\xi|^{-2 \eta}\left|\widehat{f}(\xi)\right|^2 d \xi \leq \rho^{2 \eta-2 r-n}\|f\|_{\dot{H}^{-\eta}}^2.
$$
Hence $\Theta_r(f)\lesssim\|f\|_{\dot{H}^{-\eta}}^2$ and $P_r(f)_{+}<\infty$ for $r=\eta-n/2$, and $r(f)_+\geq\eta-n/2$. If $r^*(f)$ exists, then $r^*(f)\geq\eta-n/2$.
}
\end{example}

\begin{example} \label{Example:2.3}
\rm{
We define
$$
Y^q:=\Big\{ f\in S': \widehat{f}\in L^1_{\mathrm{loc}},\quad \|f\|_{Y^q}:= \operatorname*{ess\,sup}_{\xi\in\mathbb{R}^n\setminus\{0\}} |\xi|^q|\widehat{f}(\xi)|<\infty \Big\},
$$
where $S'\left(\mathbb{R}^n\right)$ is the space of tempered distributions. If $f \in L^2 \cap Y^q$ and $q<n/2$, then
$$
\rho^{-2r-n}\int_{B(\rho)}|\widehat{f}(\xi)|^2d\xi \leq C\|f\|_{Y^q}^2\rho^{-2r-n}\int_0^\rho s^{n-1-2q}ds\leq C\|f\|_{Y^q}^2\rho^{-2q-2r}.
$$
Thus $\Theta_r(f)\lesssim\|f\|_{Y^q}^2$ and $P_r(f)_{+}<\infty$ for $r=-q$, and $r(f)_+\geq-q$. If $r^*(f)$ exists, then $r^*(f)\geq-q$.
}
\end{example}

The role of the upper decay indicators is illustrated by the following proposition.
\begin{proposition} \label{decay_1}
Let $r \in [-n/2, \infty)$ and $a(t)$ be positive, continuous and satisfy $a(t) \notin L^1([0, \infty))$. Assume that $u$ is a solution to $u_t+a(t)(-\Delta)^\sigma u=0$ with the initial datum $u(s, x)=f(s, x) \in L^2, s \geq 0$. If $P_r(f)_{+}<\infty$, then
$$
\|u(t)\|_{L^2} \lesssim \left(1+\int_s^t a(\tau) d \tau \right)^{-\frac{n+2r}{4 \sigma}}\|f\|_{r,0}.
$$
\end{proposition}

\begin{proof}
Applying the partial Fourier transform with respect to the spatial variable $x$, we obtain
$$
\widehat{u}_{t}+a(t)|\xi|^{2 \sigma}\widehat{u}=0, \quad \widehat{u}(s,\xi)=\widehat{f}(s, \xi),
$$
and hence
$$
\widehat{u}(t,\xi)=\mathrm{e}^{-|\xi|^{2\sigma}\int_s^t a(\tau)d\tau}\widehat{f}(s, \xi).
$$
Let $A_s(t):=1+\int_s^t a(\tau)d\tau$. The energy identity gives
$$
\frac{d}{dt}\|u(t)\|_{L^2}^2+2a(t)\int_{\mathbb{R}^n}|\xi|^{2\sigma}|\widehat{u}(t,\xi)|^2d\xi=0.
$$
Let $\ell_0>\max\{1,(n+2r)/(2\sigma)\}$ and define $\rho(t):=\ell_0^{1/(2\sigma)}A_s(t)^{-1/(2\sigma)}$. If $A_s(t)\leq\ell_0$, the desired estimate follows from $\|u(t)\|_{L^2}\leq\|f\|_{L^2}$. Let now $A_s(t)>\ell_0$ and choose $t_1\geq s$ such that $A_s(t_1)=\ell_0$. Then $\rho(t)\leq1$ for $t\geq t_1$, and the definition of $\Theta_r(f)$ gives
$$
\int_{B(\rho(t))}|\widehat{u}(t,\xi)|^2d\xi\leq\int_{B(\rho(t))}|\widehat{f}(s, \xi)|^2d\xi\leq\Theta_r(f)\rho(t)^{n+2r}.
$$
Therefore, Fourier splitting gives
$$
\frac{d}{dt}\|u(t)\|_{L^2}^2+2\ell_0\frac{a(t)}{A_s(t)}\|u(t)\|_{L^2}^2 \lesssim a(t)\rho(t)^{2\sigma+n+2r}\|f\|_{r,0}^2.
$$
Multiplying by $A_s(t)^{2\ell_0}$ and integrating from $t_1$ to $t$, we obtain
$$
\|u(t)\|_{L^2}^2\lesssim A_s(t)^{-\frac{n+2r}{2\sigma}}\|f\|_{r,0}^2.
$$
This completes the proof.
\end{proof}

\section{Some useful tools} \label{section_tools}

\begin{lemma}[Fractional Gagliardo-Nirenberg inequality, \cite{Hajaiej2011}]
\label{Gagliardo}
Assume that $q_1,q_2,q_3\in(1,\infty)$ and $0\leq k<r$ with $r>0$. Then, it holds
$$
\|w\|_{\dot{H}_{q_1}^k} \lesssim \|w\|_{L^{q_2}}^{1-\theta} \|w\|_{\dot{H}_{q_3}^r}^{\theta}, \text{ where }\theta\in[k/r,1] \text{ fulfilling }\frac1{q_1}-\frac{k}{n} = \frac{1-\theta}{q_2} +\theta\left(\frac1{q_3}-\frac r n\right).
$$
\end{lemma}

\begin{lemma}\label{lemma:nonlinear_weight}
Assume \textbf{(A)}-\textbf{(B)} and fix $M\geq0$ such that $\mu'(t)/\mu(t)\geq-M a(t)/A(t)$ for $t\geq0$. If $a_0>0$, then for every $0\leq\gamma< 2a_0/(2+M)$ we have $a(t)^{-2}\lesssim(1+\mathcal{V}_{a,b}(0,t))^{-\gamma}$ for $t\geq0$. If $a_0=0$, the estimate holds for $\gamma=0$.
\end{lemma}
\begin{proof}
By \textbf{(B2)}, we have
$$
\left(\frac{A^2(t)}{\mu(t)}\right)' =\left(2-\frac{A(t)\mu'(t)}{a(t)\mu(t)}\right)\frac{a^2(t)}{b(t)} \leq(2+M)\frac{a^2(t)}{b(t)}.
$$
Consequently $\frac{A^2(t)}{\mu(t)}\lesssim \left(1+\mathcal{V}_{a,b}(0,t)\right)$ and, for every $\varepsilon>0$, we obtain
$$
\frac{a^2(t)}{b(t)\left(1+\mathcal{V}_{a,b}(0,t)\right)} =\frac{A^2(t)}{\mu(t)\left(1+\mathcal{V}_{a,b}(0,t)\right)}\frac{a(t)}{A(t)} \leq(2+M+\varepsilon)\frac{a(t)}{A(t)} \quad \text{ for sufficiently large } t.
$$
Integrating gives $\left(1+\mathcal{V}_{a,b}(0,t)\right)\lesssim A(t)^{2+M+\varepsilon}$. On the other hand, \textbf{(A2)} implies $a(t)\geq a(0)A(t)^{a_0}$. Hence, $a(t)^{-2}\lesssim \left(1+\mathcal{V}_{a,b}(0,t)\right)^{-\gamma}$ whenever $\gamma(2+M+\varepsilon)<2a_0$. Since $\varepsilon>0$ is arbitrary, this proves the claimed range of $\gamma$. When $a_0=0$, we use $a(t)\geq a(0)>0$.
\end{proof}
\begin{proposition}\label{proposition:inf_beta}
Assume that $a(t),b(t)$ and $m(t)$ satisfy Assumptions \textbf{(A)}--\textbf{(B)}--\textbf{(C)}. Let $\omega\in[1,2]$, let $\beta$ be as in \eqref{eq:inf_beta}, and assume that $p>\overline p(\omega,\beta)$. Then there exist constants $\delta>0$ and $C>0$, independent of $t$, such that
$$
\frac{\chi(t)^{p-1}}{a^2(t)}\left(1+\mathcal{V}_{a,b}(0,t)\right)^{-\frac{n}{2\omega\sigma}(p-1)}\leq C\left(1+\mathcal{V}_{a,b}(0,t)\right)^{-1-\delta}, \quad t\geq0.
$$
\end{proposition}

\begin{proof}
By \textbf{(B4)}, we get $\mathcal{V}_{a,b}(0,t)\to\infty$ as $t\to\infty$. Since $m^2(t)\geq0$ and
$$
\frac{1+\mathcal{V}_{a,b}(0,t)} {\mathcal{V}_{a,b}(0,t)} \to 1,
$$
the definition of $\beta$ gives
$$
\liminf_{t\to\infty} \frac{\left(1+\mathcal{V}_{a,b}(0,t)\right)m^2(t)} {a^2(t)} =\beta.
$$
For $t>0$, we have
$$
\frac{-\log\chi(t)} {\log\left(1+\mathcal{V}_{a,b}(0,t)\right)} = \frac{ \displaystyle\int_0^t \frac{\left(1+\mathcal{V}_{a,b}(0,s)\right)m^2(s)} {a^2(s)} \frac{a^2(s)} {b(s)\left(1+\mathcal{V}_{a,b}(0,s)\right)}ds }{ \displaystyle\int_0^t \frac{a^2(s)} {b(s)\left(1+\mathcal{V}_{a,b}(0,s)\right)}ds }.
$$
The denominator tends to infinity. Integrating the eventual lower bounds supplied by the preceding liminf, and dividing by the denominator, yields
$$
\liminf_{t\to\infty} \frac{-\log\chi(t)} {\log\left(1+\mathcal{V}_{a,b}(0,t)\right)} \geq\beta,
$$
with the same interpretation when $\beta=\infty$. Moreover, Lemma \ref{lemma:nonlinear_weight} implies
$$
\liminf_{t\to\infty} \frac{2\log a(t)} {\log\left(1+\mathcal{V}_{a,b}(0,t)\right)} \geq\frac{2a_0}{2+M}.
$$
For finite $\beta$, it follows that
$$
\begin{aligned}
\liminf_{t\to\infty} \bigg\{ \frac{2\log a(t)-(p-1)\log\chi(t)} {\log\left(1+\mathcal{V}_{a,b}(0,t)\right)} +\frac{n}{2\omega\sigma}(p-1) \bigg\} &\geq \frac{2a_0}{2+M} +\left(\frac{n}{2\omega\sigma}+\beta\right)(p-1) >1,
\end{aligned}
$$
where the strict inequality follows from $p>\overline p(\omega,\beta)$. If $\beta=\infty$, the left-hand side in the previous estimate is infinite because $p>1$. Thus, in either case there exist $\delta>0$ and $t_0>0$ such that
$$
2\log a(t)-(p-1)\log\chi(t)+\frac{n}{2\omega\sigma}(p-1) \log\left(1+\mathcal{V}_{a,b}(0,t)\right) \geq (1+\delta)\log\left(1+\mathcal{V}_{a,b}(0,t)\right), \quad t\geq t_0.
$$
This proves the estimate for $t\geq t_0$. The case $0\leq t\leq t_0$ follows by continuity.
\end{proof}

\end{document}